\documentclass[12pt]{article}
\usepackage{geometry}                
\usepackage{graphicx}
\usepackage{amssymb}
\usepackage{amsmath}
\usepackage{epstopdf}
\usepackage{hyperref}
\usepackage{nomencl}
\usepackage{trfsigns}
\usepackage{arydshln}
\usepackage{amsthm}
\usepackage{subfig}
\usepackage{color}

\usepackage{float}
\usepackage{placeins}

\def\R2Lurl#1#2{\mbox{\href{#1}{\tt #2}}}

\title{An Exact SIR Series Solution and an Exploration of the Related Parameter Space}
\author{Daniel P. Hobbs\\
{\normalsize\itshape{School of Mathematics and Statistics}}\\ 
{\normalsize\itshape{Rochester Institute of Technology, Rochester, NY 14623, USA}} \\ 
}
\date{}

\begin{document}
\maketitle
\begin{abstract}
\baselineskip = 2.0\baselineskip
  A convergent power series solution is obtained for the SIR model, using an asymptotically motivated gauge function. For certain choices of model parameter values, the series converges over the full physical domain (i.e., for all positive time). Furthermore, the radius of convergence as a function of nondimensionalized initial susceptible and infected populations is obtained via a numerical root test.
\end{abstract}
\baselineskip = 2.0\baselineskip

\maketitle

\section{Introduction and Background}
\subsection{Defining the Problem}
The SIR model is a classic compartmental model in mathematical epidemiology in which we have the time-evolving variables $S$ susceptible population, $I$ infected population, and $R$ recovered population. It was originally proposed by Kermack and McKendrick~\cite{Kermack}. Before explaining the motivation of the form of the equations, it is noted that the SIR equations are the following.

\begin{subequations}
\begin{equation}
    \frac{dS}{dt}\boldsymbol{=} -rSI,\,\,\,\,\,\,\, S(0)\boldsymbol{=} S_{0},\,\,\,\,\,\,\, 0\leq t\leq\infty
\end{equation}

\begin{equation}
    \frac{dI}{dt}\boldsymbol{=} rSI-\alpha I,\,\,\,\,\,\,\, I(0)\boldsymbol{=} I_{0}
\end{equation}

\begin{equation}
    \frac{dR}{dt}\boldsymbol{=}\alpha I,\,\,\,\,\,\,\, R(0)\boldsymbol{=} 0
\end{equation}

\label{SIReqn}
\end{subequations}

In general, it is assumed that the time derivative $dS/dt$ is negatively proportional to $S$ and $I$ since the susceptible population is assumed to decline at a greater rate with a larger infected population while its rate of decline should be proportional to the susceptible population itself. More specifically, one would expect this the time derivative to be negatively proportional to the multiplicative product of $S$ and $I$ since this is in a sense a direct measure of the interaction between the susceptible and infected populations. Furthermore, for every amount contributing to the decline of the susceptible population, the same amount $rSI$ should positively contribute to the growth $dI/dt$ of the infected population.

However, it is also possible for some of the infected to recover, and so, it is assumed there is a linear subtraction $\alpha I$ from the infected rate of increase. Note that this linear quantity proportional to the infected ultimately gets integrated to give the recovered population. In summary, this overall description leads us to the classic SIR equations (\ref{SIReqn}).

\subsection{Analytical Solution Techniques}

In general, it is possible to obtain an analytical power series solution to a system of nonlinear ordinary differential equation such as the SIR model. Even as a power series, if implemented by a clever strategy to give a useful solution, then it is reasonably possible that this will lead to greater computational efficiency and accuracy compared to the usual numerical methods such as the Runge-Kutta method (often abreviated as RK4). Since this means avoiding the necessity in numerical step methods of having to use many small time steps, the series approach can make an enormous difference if it is being implemented as part of a much larger problem utilizing a Monte Carlo simulation for example. However, the strategy in question must deal with the issue of a time power series usually having a radius of convergence $\rho$ in which the power series is exact for $|t|<\rho$ even though the physical domain is $0\le t \le\infty$. Sometimes, such as in the case of the SIR model, this problem can be solved through approximant techniques which is what was previously done by Barlow and Weinstein~\cite{BarlowWeinstein}, and the goal of this paper is to accomplish the same thing through exact techniques.
\subsection{Goal of the Project}

From Barlow and Weinstein~\cite{BarlowWeinstein}, an autonomous equation of only $S$ the susceptible population is obtained. This was obtained by considering the ratio of the infected and susceptible population rates and integrating to obtain an explicit expression of the infected population in terms of the susceptible population.
\begin{equation}
    I\boldsymbol{=}\frac{\alpha}{r}\ln\left(\frac{S}{S_{0}}\right)-S+S_{0}+I_{0}
\end{equation}

Upon substitution back into the usual susceptible population derivative equation, we now have the following first order ODE.
\begin{equation}
    \frac{dS}{dt}\boldsymbol{=}\beta S+r S^{2} - \alpha S\ln(S),\,\,\,\,\, \beta\boldsymbol{=} \alpha\ln(S_{0})-r(S_{0}+I_{0}),\,\,\,\,\, S(0)\boldsymbol{=} S_{0},\,\,\,\,\, 0\leq t\leq\infty
    \label{SWB}
\end{equation}

Since the total sum of the susceptible, infected, and recovered populations is a conserved quantity, the following equation also follows with the assumption that the initial recovered population is $R_{0}\boldsymbol{=} 0$.

\begin{equation}
    R\boldsymbol{=} S_{0}+I_{0}-S-I\boldsymbol{=}-\frac{\alpha}{r}\ln\left(\frac{S}{S_{0}}\right)
    \label{4}
\end{equation}

As a result, our only dynamical unknown is the susceptible population $S$, and therefore, its solution will fully determine the solutions of the infected and recovered populations $I$ and $R$. So, this is precisely why this paper will focus exclusively on the exact analytical solution of $S$.

To solve the SIR model, we want to construct a power series for $S$ in a similar manner to what was done by Naghshineh et al~\cite{Nastaran}. To explain further, there are a number of instances in which if an appropriate independent variable substitution is utilized in an ordinary differential equation, then the new resulting series can sometimes converge over the full physical domain, leading to an exact analytical solution. In fact, there is precedence for this, and for example, this is exactly what was accomplished in two problems seen in the previously mentioned paper by Naghshineh et al~\cite{Nastaran}.

Building off the dominant balance argument of~\cite{BarlowWeinstein} leads to the exponential pattern describing a nondimensionalized quantity directly through a one-to-one function of $S$ giving... $A_{0}+A_{1}e^{\lambda T}+A_{2}e^{2\lambda T}+...$ with $T$ being nondimensionalized time, and so, this motivates us to use an exponential type gauge substitution. We will see later that a modification to this gauge will be required in order for terms to be determined exactly by initial conditions. Unlike in Nagshineh et al, we have four parameters which is hardly ideal, and therefore, we need to determine out how and if we can collapse these parameters in order to more easily constrain the eventual convergence survey of our improved series. In particular, one of the main goals of this paper is to determine the parameter space for which a convergent exact solution is obtainable.

The report is organized as follows. In section 2, referencing a simplification of the problem through a collapse of parameters (which is explicitly explained in appendices A and B), it is shown that a calculation of the direct time series is divergent within the physical domain. In section 3, an asymptotic analysis is performed motivating the exponential gauge, and it is shown why it is preferred to work with a shifted version of this series. In section 4, the results are presented which includes contour plots illustrating the parameter space, and in section 5, there is a discussion involving some asymptotic analysis along with ideas for future work. In section 6, the conclusions are presented.

\section{Problem Formulation and Divergent Series}

As a brief summary of key results which are derived in the appendices (see Appendices A and B), the problem is structured as follows with an associated collapse of parameters. Also, note that as expected, equation (\ref{xiTL}) produces a monotonically decreasing function (see Appendix D).
\begin{subequations}
\begin{equation}
        \frac{d\xi}{dT}\boldsymbol{=} L\xi^{2}-\xi\ln(\xi),\,\,\,\,\,\,\, \xi(0)\boldsymbol{=}\xi_{0}\boldsymbol{=}e^{\tilde{S_{0}}+\tilde{I_{0}}},\,\,\,\,\,\,\, 1\leq\xi\leq\infty
        \label{xiTL}
    \end{equation}
    \begin{equation}
        T\boldsymbol{=}\alpha t,\,\,\,\,\,\,\, \xi\boldsymbol{=}\frac{\tilde{S}}{L},\,\,\,\,\,\,\, L\boldsymbol{=}\tilde{S_{0}}e^{-\tilde{S_{0}}-\tilde{I_{0}}}
    \end{equation}
\begin{equation}
        \tilde{S}\boldsymbol{=}\frac{rS}{\alpha},\,\,\,\,\,\,\, \tilde{I}\boldsymbol{=}\frac{rI}{\alpha}
    \end{equation}
\end{subequations}

For convenience, the substitution $\xi\boldsymbol{=}e^{V}$ (which was motivated for a special case in Appendix E) will be used for the general case, and this naturally assumes the initial condition $V_{0}\boldsymbol{=}\ln(\xi_{0})\boldsymbol{=}\tilde{S_{0}}+\tilde{I_{0}}$. Incidentally, the connections between the main equations can be summarized in the following flowchart.
\begin{equation}
    \frac{dV}{dT}\boldsymbol{=} L e^{V}-V,\,\,\,\,\,\,\, V(0)\boldsymbol{=} V_{0}\boldsymbol{=}\tilde{S_{0}}+\tilde{I_{0}}
    \label{VEQN}
\end{equation}

\begin{figure}[h!]
    \centering
    \includegraphics[width=0.8\linewidth]{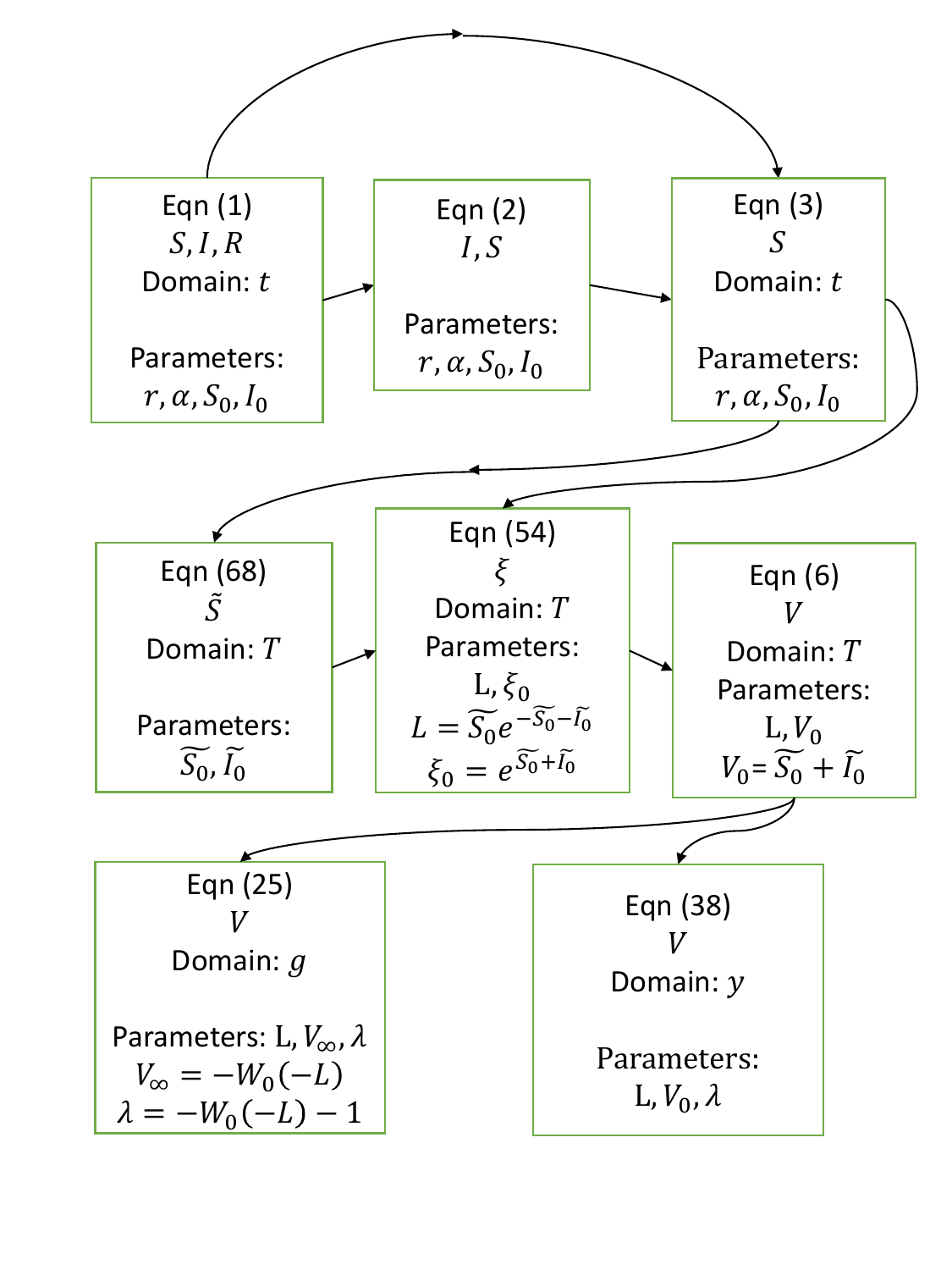}
    \caption{This intent of this flowchart is to briefly illustrate how the various equations of dynamical variables and their parameters connect to one another.}
    \label{fig:enter-label}
\end{figure}

The original motivation for using the substitution $\xi\boldsymbol{=}e^{V}$ was that it leads directly to an ODE (\ref{VEQN}) which is simpler algorithmically to solve. Having said that, it turns out that $V$ does have a "physical" meaning, and in general, it is always the sum of the nondimensionalized susceptible and infected populations. This can be shown by substituting its definition into equations (\ref{SWB}), and (\ref{4}) while referencing equations (5b) and (5c).

\begin{equation}
    V\boldsymbol{=}\tilde{S}+\tilde{I}\boldsymbol{=}\tilde{S_{0}}+\tilde{I_{0}}-\tilde{R},\,\,\,\,\,\,\, \tilde{R}\boldsymbol{=}\frac{rR}{\alpha}
\end{equation}

Next, we want to focus on how to obtain solutions to $V$. To reiterate, substituting $\xi\boldsymbol{=}e^{V}$ in equation (5a) results in an equation that is simpler to solve as a series solution in terms of the related algorithm. Implementing the series algorithm, we will assume the following series setup below. Furthermore, solving for the related coefficients gives us the following recursion formulas. Note that the $B_{n}$ coefficients always contain a single factor $\xi_{0}$. Depending on the exact values of $\tilde{S_{0}}$ and $\tilde{I_{0}}$, $\xi_{0}$ will have a tendency to be numerically large. Consequently, it might be preferable to instead work with the $C_{n}$ coefficients without referencing the $B_{n}$. This will perhaps result in less numerical error when $\xi_{0}$ is large.

\begin{equation}
    C_{n}\boldsymbol{=} L B_{n}
\end{equation}

\begin{equation}
    A_{0}\boldsymbol{=}\tilde{S_{0}}+\tilde{I_{0}}
\end{equation}

\begin{equation}
    B_{0}\boldsymbol{=}e^{\tilde{S_{0}}+\tilde{I_{0}}}\boldsymbol{=}\xi_{0}
\end{equation}

\begin{equation}
    C_{0}\boldsymbol{=}L B_{0}\boldsymbol{=}\tilde{S_{0}}
\end{equation}

\begin{equation}
    A_{n+1}\boldsymbol{=}\frac{L B_{n}-A_{n}}{n+1}\boldsymbol{=}\frac{C_{n}-A_{n}}{n+1}
\end{equation}

\begin{equation}
    B_{n+1}\boldsymbol{=}\frac{1}{n+1}\sum_{j\boldsymbol{=} 0}^{n} (j+1)A_{j+1}B_{n-j}
\end{equation}

\begin{equation}
    C_{n+1}\boldsymbol{=}\frac{1}{n+1}\sum_{j\boldsymbol{=} 0}^{n} (j+1)A_{j+1}C_{n-j}
\end{equation}

\begin{equation}
    V\boldsymbol{=}\sum_{n\boldsymbol{=} 0}^{\infty} A_{n} T^{n}
    \label{8}
\end{equation}

\begin{equation}
    \xi\boldsymbol{=} e^{V}\boldsymbol{=} \sum_{n\boldsymbol{=} 0}^{\infty} B_{n} T^{n}
\end{equation}

\begin{equation}
    \tilde{S}\boldsymbol{=} L\xi\boldsymbol{=} \sum_{n\boldsymbol{=} 0}^{\infty} C_{n} T^{n}
\end{equation}

Note that equations (13) and (14) can be found from a previously published result~\cite{Gibbons}. This result of finding the power series of the exponential of a power series will be used in later in the paper as well. Also, the first few coefficients are listed in Appendix F. Finally, note that a radius of convergence can be calculated with the root test.

\begin{equation}
    \rho\boldsymbol{=}\lim_{n \to \infty} |A_{n}|^{-1/n}
\end{equation}

\begin{figure}[h!]
    \centering
    \includegraphics[width=1\linewidth]{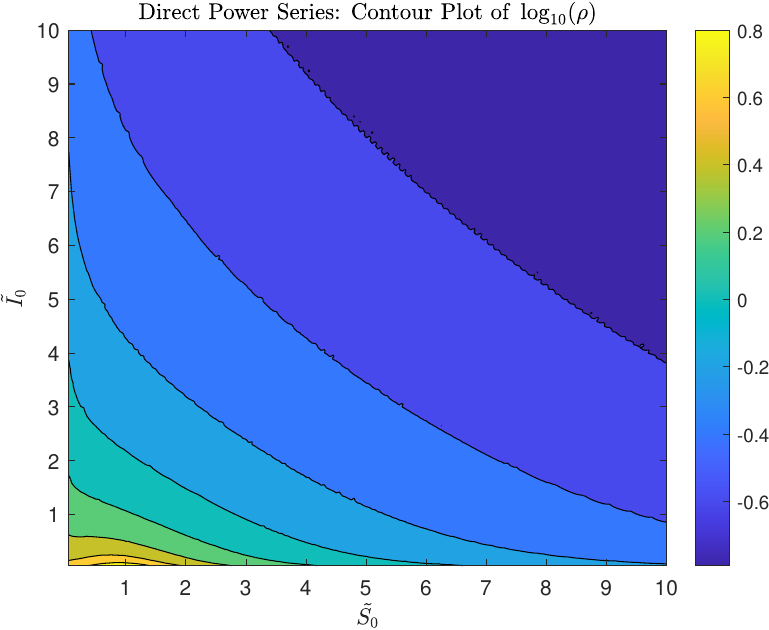}
    \caption{For the direct time series contour plot of $\log_{10}(\rho)$ of the power series for $V$ given by (8), $N\boldsymbol{=} 300$ was used. In other words, $\rho\approx A_{300}^{-1/300}$ was used.}
    \label{fig: Figure 1}
\end{figure}

It can be observed that this always results in a finite radius of convergence $\rho$ which ends up being a function of $\tilde{S_{0}}$ and $\tilde{I_{0}}$ in this system. In short, there is always a divergence within the physical domain. This is exactly what we see in this case (Figure \ref{fig: Figure 1}) which shows the numerically calculated $\rho$. Possibly, arbitrarily large radii can be found when $\tilde{I_{0}}$ is very small while $\tilde{S_{0}}=1$, but most likely, a truly infinite radius of convergence will never occur even under this very restricted set of circumstances.

\begin{figure}[h!]
    \centering
    \includegraphics[width=1\linewidth]{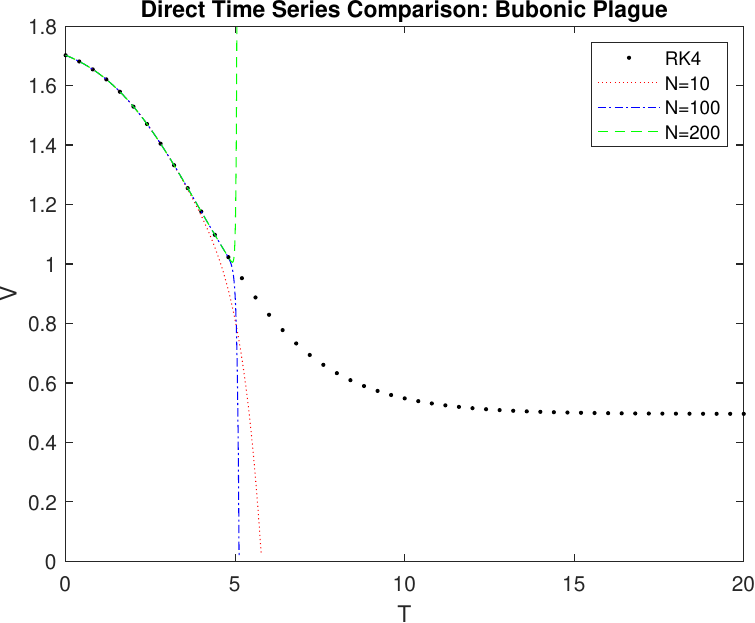}
    \caption{This is a direct nondimensionalized time series comparison for the 1966 bubonic plague which took place in Eyam, England with the parameters $S_{0}\boldsymbol{=} 254$, $I_{0}\boldsymbol{=} 7$, $r\boldsymbol{=} 0.0178$, and $\alpha\boldsymbol{=} 2.73$~\cite{BarlowWeinstein, Khan}. This can be rewritten in the nondimensionalized quantities $(\tilde{S_{0}}, \tilde{I_{0}}) \boldsymbol{=} (1.656117, 0.045641)$.}
    \label{fig: Figure 2}
\end{figure}

For further illustration, the 1966 bubonic plague which can be identified on Figure \ref{fig: Figure 1} as $(\tilde{S_{0}},\, \tilde{I_{0}})\boldsymbol{=} (1.656117,\, 0.045641)$, when plotted as a nondimensionalized time series in Figure \ref{fig: Figure 2} shows divergent behavior indicating a radius of convergence roughly equal to $5$. Any possible series approach to solving the SIR will have to contend with circumventing this divergent behavior. In fact, our goal is to do precisely this by constructing a resummation of (\ref{8}) which will consequently converge over the full physical domain.

\newpage
\section{Convergent Resummation and Asymptotics}
\subsection{Dominant Balance of V}

In search of a natural gauge variable, it seems best to do a dominant balance analysis. So, near the limit of $T\rightarrow\infty$, we will assume the function obeying equation (\ref{VEQN}), is approximately a combination of $V_{\infty}\boldsymbol{=}\ln(\xi_{\infty})$ and a function $G(T)$ which is much smaller than $V_{\infty}$. This argument can be made exact in some sense too by including higher order powers of $G$ in the combination.
\begin{equation}
    V\approx V_{\infty}+G\boldsymbol{=}\ln(\xi_{\infty})+G
\end{equation}

Substituting this into equation (\ref{VEQN}) now gives us the following.
\begin{equation}
    \frac{dG}{dT}\boldsymbol{=} L e^{(V_{\infty}+G)}-V_{\infty}-G\boldsymbol{=} L\xi_{\infty}e^{G}-\ln(\xi_{\infty})-G
    \label{Vg1}
\end{equation}

Note that $e^x\boldsymbol{=} 1+x+x^{2}/2+x^{3}/6+...$, and so, the exponential contribution gives the following.
\begin{equation*}
    L\xi_{\infty}e^{G}-\ln(\xi_{\infty})-G\boldsymbol{=} L\xi_{\infty}(1+G+G^{2}/2+G^{3}/6+...)-\ln(\xi_{\infty})-G
\end{equation*}
\begin{equation*}
    \boldsymbol{=} L\xi_{\infty}-\ln(\xi_{\infty})+(L\xi_{\infty}-1)G+L\xi_{\infty}(G^{2}/2+G^{3}/6+...)
\end{equation*}

Keeping only first order terms, we get the following.
\begin{equation}
    L\xi_{\infty}e^{G}-\ln(\xi_{\infty})-G\boldsymbol{=} L\xi_{\infty}-\ln(\xi_{\infty})+(L\xi_{\infty}-1)G
\end{equation}

From equation (\ref{xiInfinityEQN}) in Appendix C, we know that $L\xi_{\infty}-\ln(\xi_{\infty})\boldsymbol{=} 0$. This occurs since this expression is $\frac{1}{L}\frac{d\xi}{dT}$, and as time approaches infinity, the derivative will approach zero. So, by simple substitution into equation (\ref{Vg1}), we obtain the following.
\begin{equation}
    \frac{dG}{dT}\boldsymbol{=} (L\xi_{\infty}-1)G
    \label{Vg2}
\end{equation}

Consequently, equation (\ref{Vg2}) results in the dominant balance giving us the simplest nontrivial solution below which is valid near $T\rightarrow\infty$. Note that $G$ and $g$ are different in order to distinguish between having the proportionality constant $a_{1}$ versus dropping it.
\begin{equation}
    G\boldsymbol{=} a_{1} g, \,\,\,\,\,\,\, g\boldsymbol{=} e^{\lambda T}
    \label{Vg3}
\end{equation}

Incidentally, the constant $\lambda$ was solved in terms of the Lambert function $W_{0}(x)$ which mathematically is the solution to the generic problem $f e^f\boldsymbol{=} x$.
\begin{equation}
    \lambda\boldsymbol{=} L\xi_{\infty}-1\boldsymbol{=} -W_{0}(-L)-1
\end{equation}

Note that a dominant balance analysis of equation (\ref{xiTL}) produces the exact same result (see Appendix G). If you already know $a_{1}$, it is possible to motivate higher order correction terms $a_{2}g^{2}$, $a_{3}g^{3}$, .... by the same procedure as before. In other words, these terms can get added into equation (19) to make it exact and then once again substituted into equation (\ref{VEQN}). This is what motivates the gauge $g\boldsymbol{=} e^{\lambda T}$ that, when used in a power series of terms $g^{n}$ will naturally represent the asymptotic ordering of the dominant balance.

\subsection{Straight Exponential Gauge Variable Substitution}

Note that equation (\ref{Vg3}) transforms the domain by mapping $T\boldsymbol{=}\infty$ to $g\boldsymbol{=} 0$ and $T\boldsymbol{=} 0$ to $g\boldsymbol{=} 1$. So, if a power series solution were constructed from powers of $g$, then it would be sufficient to prove a fully convergent solution if the radius of convergence is greater than or equal to one. Also, note it will transform equation (\ref{VEQN}) into the following.

\begin{equation}
    \frac{dV}{dg}\boldsymbol{=}\frac{L e^{V}-V}{\lambda g},\,\,\,\,\,\,\, V(g\boldsymbol{=}0)\boldsymbol{=} V_{\infty}\boldsymbol{=}\ln(\xi_{\infty})\boldsymbol{=}-W_{0}(-L),\,\,\,\,\,\,\, 0\leq g\leq 1
\end{equation}

Assuming now that $V$ is a power series of $g$, it follows with the initial condition that we can write the power series as the following.
\begin{equation}
    V\boldsymbol{=}\sum_{n\boldsymbol{=}0}^{\infty}a_{n}g^{n},\,\,\,\,\,\,\,\,\,\,\,\, V(g\boldsymbol{=}0)\boldsymbol{=}a_{0}\boldsymbol{=} V_{\infty}\boldsymbol{=}\ln(\xi_{\infty})\boldsymbol{=} -W_{0}(-L)
\end{equation}

\begin{equation}
    \tilde{S}\boldsymbol{=} L\xi\boldsymbol{=} L e^V\boldsymbol{=}\sum_{n\boldsymbol{=} 0}^{\infty}b_{n}g^{n}
\end{equation}

\begin{equation}
    b_{0}\boldsymbol{=} L e^{a_{0}}\boldsymbol{=} L\xi_{\infty}\boldsymbol{=}\tilde{S}_{\infty}\boldsymbol{=}-W_{0}(-L)
\end{equation}
\begin{equation}
    a_{n+1}\boldsymbol{=}\frac{b_{n+1}-a_{n+1}}{(n+1)\lambda}
\end{equation}
\begin{equation}
    a_{n+1}\boldsymbol{=}\frac{b_{n+1}}{(n+1)\lambda+1}
    \label{AgRecursive}
\end{equation}
\begin{equation}
    b_{n+1}\boldsymbol{=}\frac{1}{n+1}\sum_{j\boldsymbol{=} 0}^{n} (j+1)a_{j+1}b_{n-j}
    \label{BgRecursive}
\end{equation}

Note that the previous recursive equations make things a little awkward computationally compared to what was done before because it will not be possible to recursively solve for $a_{1}$ in a direct manner. To illustrate this, let us begin to carry out the recursion. Doing the first iteration of equations (\ref{AgRecursive}) and (\ref{BgRecursive}) gives the following.

\begin{equation}
    a_{1}\boldsymbol{=}\frac{L b_{1}}{\lambda+1}\boldsymbol{=}\frac{L b_{1}}{L\xi_{\infty}}\boldsymbol{=}\frac{b_{1}}{\xi_{\infty}}
\end{equation}
\begin{equation}
    b_{1}\boldsymbol{=}a_{1}b_{0}\boldsymbol{=}\xi_{\infty}a_{1}
\end{equation}

As we can see, in the first iteration, both recursions ended up telling us equivalent information dependent on $A_{1}$ without us knowing what it is. This issue can only be resolved by noting that at $g\boldsymbol{=} 1$, we arrive at $T\boldsymbol{=} 0$ meaning $V=\ln(\xi_{0})\boldsymbol{=}\tilde{S_{0}}+\tilde{I_{0}}$. However, this will require an infinite number of coefficients to be simultaneously considered giving us the following.

\begin{equation}
    a_{0}+a_{1}+a_{2}+a_{3}+...\boldsymbol{=}\ln(\xi_{0})\boldsymbol{=}\tilde{S_{0}}+\tilde{I_{0}}
    \label{AgFinal}
\end{equation}

In order to solve this for $a_{1}$ and other coefficients, it must be solved iteratively through substitution in terms of $a_{1}$ for each coefficient $a_{n>1}$ from the iterative equations (\ref{AgRecursive}) and (\ref{BgRecursive}). Finally, upon substitution into equation (\ref{AgFinal}), one would obtain a power series equation of $a_{1}$ that would only be solvable through a root finding method. Technically, this requires an infinite number of substitutions, but we can always do the root finding with a cut off at some power of $a_{1}$.

This is already computationally annoying, but it is nonetheless possible to implement. To make this procedure more clear, the following could be defined.

\begin{equation}
    E_{n}\boldsymbol{=}\frac{a_{n}}{a_{1}^{n}}
\end{equation}

\begin{equation}
    F_{n}\boldsymbol{=}\frac{b_{n}}{a_{1}^{n}}
\end{equation}

This now makes it possible to solve for the $E_{n}$ and $F_{n}$ coefficients iteratively (see Appendix H). Consequently, it also becomes possible to solve for $a_{1}$ with (\ref{AgFinal}) transforming into a power series equation of $a_{1}$. However, there are computational complications that result in less than ideal numerical accuracy in calculated values of $V(g)$. Consequently, we want to modify our gauge variable in some way to remove this difficulty.

\subsection{Solution Written as a Power Series of Shifted Gauge Variable}

Adding a constant to the established dominant balance, if we choose our criteria to map $T\boldsymbol{=} 0$ to $y\boldsymbol{=} 0$ and $T\boldsymbol{=}\infty$ to $y\boldsymbol{=} 1$, then the following must be our gauge variable.

\begin{equation}
    y\boldsymbol{=} 1-g\boldsymbol{=} 1-e^{\lambda T}
\end{equation}

This has the major advantage of restoring our original initial condition. As a result, we can rewrite equation (\ref{VEQN}) as the following.

\begin{equation}
    \frac{dV}{dy}\boldsymbol{=}\frac{V-L e^{V}}{\lambda(1-y)},\,\,\,\,\,\,\,\,\,\,\, V(y\boldsymbol{=} 0)\boldsymbol{=} V_{0}\boldsymbol{=}\tilde{S_{0}}+\tilde{I_{0}},\,\,\,\,\,\,\,\,\,\,\,\, 0\leq y\leq 1
    \label{VEQNy}
\end{equation}

Once again assuming a power series expansion for $V$ with newly defined coefficients $A_{n}$, these can be solved for along with the related coefficients thereby giving us the following recursion formulas while once again making $C_{n}\boldsymbol{=} L B_{n}$ an option to work with.

\begin{equation}
    A_{0}\boldsymbol{=} V_{0}
\end{equation}

\begin{equation}
    B_{0}\boldsymbol{=} e^{V_{0}}\boldsymbol{=}e^{\tilde{S_{0}}+\tilde{I_{0}}}\boldsymbol{=}\xi_{0}
\end{equation}

\begin{equation}
    C_{0}\boldsymbol{=}\tilde{S_{0}}
\end{equation}

\begin{equation}
    A_{n+1}\boldsymbol{=}\frac{(n\lambda+1)A_{n}-L B_{n}}{(n+1)\lambda}\boldsymbol{=}\frac{(n\lambda+1)A_{n}-C_{n}}{(n+1)\lambda}
\end{equation}

\begin{equation}
    B_{n+1}\boldsymbol{=}\frac{1}{n+1}\sum_{j\boldsymbol{=} 0}^{n} (j+1)A_{j+1}B_{n-j}
\end{equation}

\begin{equation}
    C_{n+1}\boldsymbol{=}\frac{1}{n+1}\sum_{j\boldsymbol{=} 0}^{n} (j+1)A_{j+1}C_{n-j}
\end{equation}

In contrast to the approximant approach, our recursion is completely defined. As aside, note that it will be preferable to work with the iterative combination of the $A_{n}$ and $C_{n}$ coefficients. The first few coefficients are calculated in Appendix I.

\section{Results}
\subsection{Radius of Convergence and Error Contour Plots}
Let us now reconsider the 1966 bubonic plague, but now, we will do so with the shifted exponential series resummation as seen in Figure 4. Note that with $N=35$, the new series prediction is indistinguishable from the original RK4 prediction. Hence, for this specific case, we have bypassed the original radius of convergence that existed in its original direct time series form. So, the next logical question is this. In the context of the resummation, what other cases also achieve full physical convergence?
\begin{figure}[h!]
    \centering
    \includegraphics[width=1\linewidth]{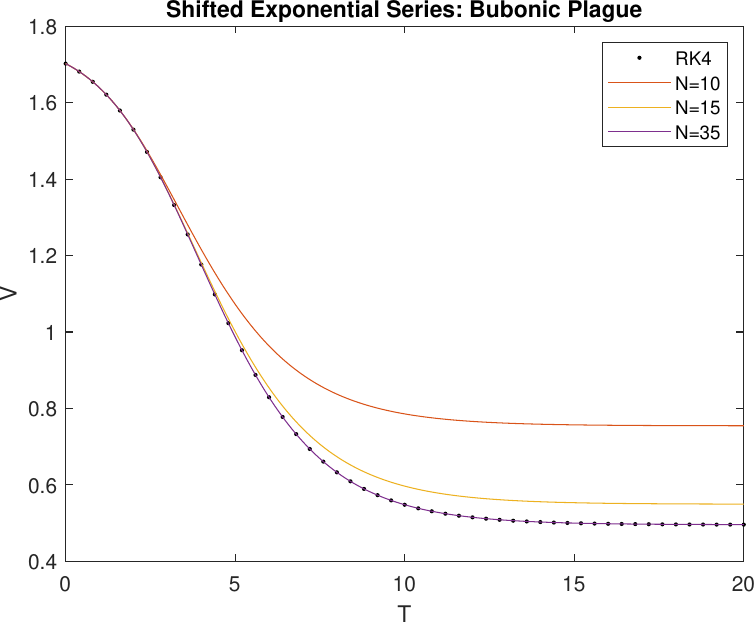}
    \caption{The above plots are the series predictions of V for the 1966 Bubonic Plague with $(\tilde{S_{0}}, \tilde{I_{0}}) \boldsymbol{=} (1.656117, 0.045641)$, $N\boldsymbol{=} 10$, $N\boldsymbol{=} 15$, and $N\boldsymbol{=} 35$.}
    \label{fig: Figure 3}
\end{figure}

To answer this question, we want to systematically investigate the radius of convergence of equation (\ref{VEQNy}), and so, let us consider the parameter space $(\tilde{S_{0}},\tilde{I_{0}})$. Next, we want to see how the radius of convergence depends on the parameter space coordinates on a contour plot. To be clear, recall that the semi-infinite physical domain $0\leq T\leq\infty$ has been mapped to a new domain $0\leq y\leq 1$, and so, we only need to have our radius of convergence be greater than or equal to $1$ in order to have full convergence. In fact, this "image" will vary with the number of terms included to determine the radius of convergence, but as expected, it stabilizes for increasing $N$.

\begin{equation}
    \rho_{N}\boldsymbol{=} |A_{N}|^{-1/N}
\end{equation}

Sometime before reaching $N\boldsymbol{=} 1000$, the yellow shape of the fully convergent region is stabilized as far as the overall details look on a macro level. In general, we end up with a "hershey kiss" shaped region (proposed by Nathan Barlow as the name of the fully convergent area). Inside the Hershey Kiss (as seen in Figure 5), the radius of convergence from the root test is determined to be numerically slightly greater than $1$.
\begin{figure}[h!]
    \centering
    \includegraphics[width=1\linewidth]{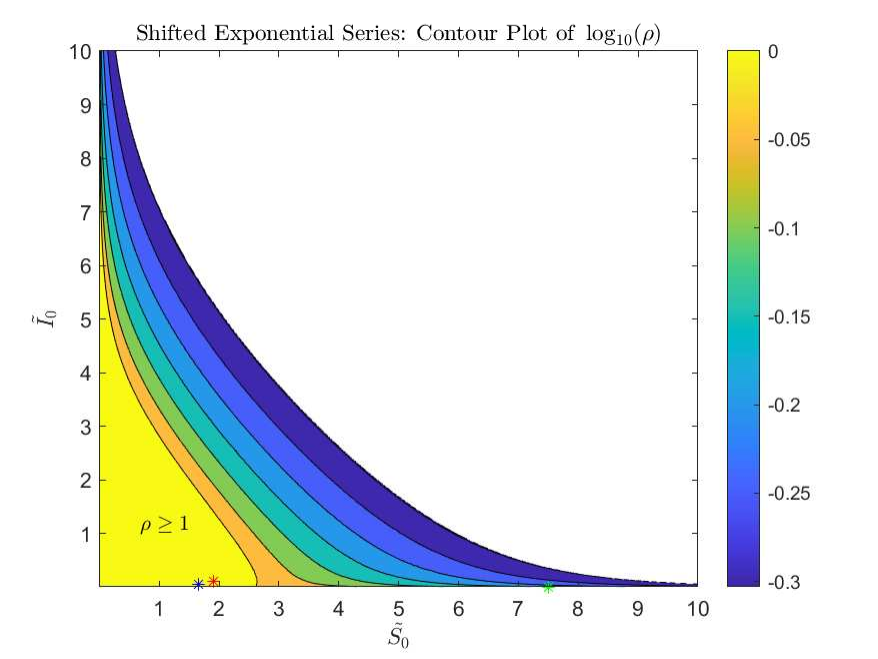}
    \caption{In the $\log_{10}(\rho)$ contour plot, $N\boldsymbol{=} 1000$ is used, and the fully convergent region (with $\rho\geq 1$) called the Hershey-Kiss region is the area shown in yellow. The blue asterisk is the Bubonic Plague at $(1.656117, 0.045641)$ with $\rho_{bubonic}\approx 1.048927$. With $r\boldsymbol{=} 0.2$, $\alpha\boldsymbol{=} 0.1$, $S_{0}\boldsymbol{=} 0.95$, and $I_{0}\boldsymbol{=} 0.05$~\cite{Rachah}, the red asterisk is Ebola at $(1.9, 0.1)$ with $\rho_{ebola}\approx 1.04929$. With $r\boldsymbol{=} 2.9236\times 10^{-5}$, $\alpha\boldsymbol{=} 0.0164$, $S_{0}\boldsymbol{=} 4206$, and $I_{0}\boldsymbol{=} 2$~\cite{BarlowWeinstein,covid}, Covid-19 in Japan is represented by the green asterisk located at $(7.497964, 0.003565)$ with $\rho_{covid}\approx 0.735302$. Note that since the basic reproductive number $r/alpha$ factors out of the analysis, it is not playing any role in the formation of the Herskey-Kiss region.}
    \label{fig: Figure 4}
\end{figure}

In summary, the shifted exponential series converges over the entire physical domain for all parameter combinations represented by the the HK (Hershey-Kiss) region. Examples such as the Ebola (red asterisk in Figure 5) and the Bubonic plague (blue asterisk in Figure 5) models are represented by points within this region.

It is also worth considering the points $(1, e-1)$ and $(e, 1)$ which both have the same corresponding values of $L\boldsymbol{=} e^{-e}$. The first point is inside the HK region whereas the second point is outside. From the root test used in Figure 5, their radii of convergence are $\rho_{1}\approx 1.049995888173082$ and $\rho_{2}\approx 0.859278414606642$ respectively.

Furthermore, note the parameter space excludes all points that are exactly on the $\tilde{S_{0}}$ axis or the $\tilde{I_{0}}$ axis. In other words, trivial cases in which $\tilde{I_{0}}\boldsymbol{=} 0$ or $\tilde{S_{0}}\boldsymbol{=} 0$ are excluded. In fact, the behavior off of these axis's is fundamentally different in some ways compared to being exactly on the axes. This is a result of any curve with a very small $\tilde{I_{0}}<<1$ (arbitrarily close to $0$) and nonzero $\tilde{S_{0}}$ will make a curved trajectory going from right to left with $\tilde{I}$ being a function of $\tilde{S}$ in which the infected eventually reaches a maximum and then begins decreasing back to $0$. This behavior will never be seen if $\tilde{I_{0}}\boldsymbol{=} 0$ which is why the on-axis behavior is different.

It is interesting to wonder why we get a Hershey Kiss curve which has the appearance of $\tilde{S}_{0,max}$ being a Gaussian curve of $\tilde{I_{0}}$, and it does seem strange that in the lower right corner, there appears to be a maximum value of $\tilde{S_{0}}$ located at $(\tilde{S_{0}}, 0)$ for which the behavior is fully convergent in the limit of $\tilde{I_{0}}\rightarrow 0$. In the context of several asymptotic approaches (see Appendices J, K, and L), this issue is further addressed later in the discussion along with appendices J and M.

Finally, as a consistency check, the common logarithm of the maximum calculated error of $V$ compared to the RK4 method produces a similar picture with the same Hershey Kiss area showing up (except now in dark blue in Figure 6). In terms of technical specifics, this implementation of RK4 used a domain of $0\leq T\leq 20$ with nondimensionalized time steps of $\Delta T\boldsymbol{=} 10^{-4}$. Doing RK4 with the nondimensionalized time domain resulted in better computational accuracy compared to implementation in the $y$ domain. As for the $y$ series approach, each series was calculated up to the $A_{1000}y^{1000}$ term. Finally, taking the absolute difference of the RK4 and the series to represent the error, the maximum error was found from within the range of values the error takes on from $0\leq T\leq 20$. Then, the pixel representing the given $(\tilde{S_{0}}, \tilde{I_{0}})$ gets color coded according to the maximum error found from its resulting scenario.
\begin{figure}[h!]
    \centering
    \includegraphics[width=1\linewidth]{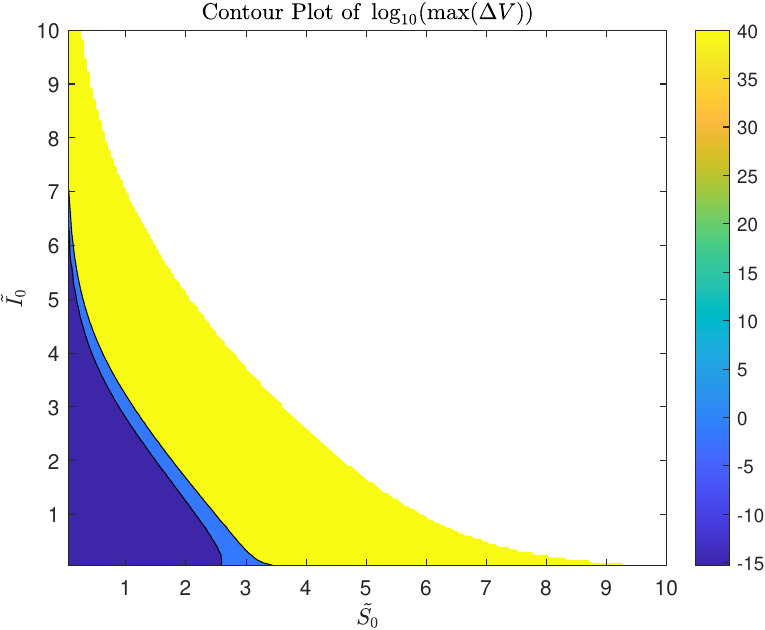}
    \caption{Each point is colored according to common logarithm of the maximum value of the error $\Delta V\boldsymbol{=} |V_{RK4}-V|$ that occurs (from the corresponding initial conditions) over the interval $0\leq y\leq 1$.}
    \label{fig: Figure 5}
\end{figure}

\subsection{Specific Error Plots of Known Disease Scenarios}

As previously noted, the initial conditions of the Bubonic Plague and Ebola are both inside the parameter space of the HK region as denoted by the blue and red asterisks respectively. However, one noted scenario of Covid-19 (in Japan) denoted by the green asterisk lies outside this convergent region, but nonetheless, our methodology can still be used to make predictions for it up until its corresponding radius of convergence. Let us first consider 3 error plots of the Bubonic plague for different values of $N$ as seen in Figure 7. 
\begin{figure}[h!]
    \centering
    \includegraphics[width=1\linewidth]{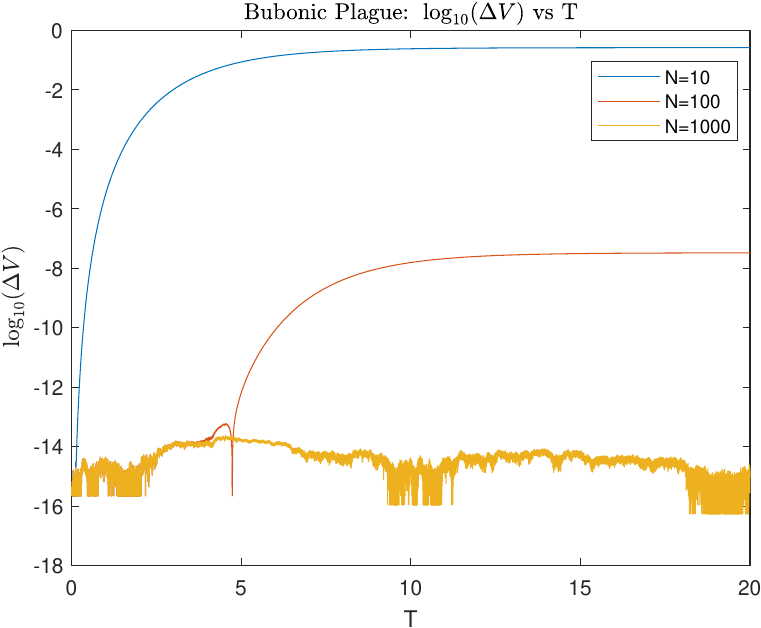}
    \caption{In this Bubonic Plague plot, we have the logs of the errors of V calculated (in comparison to RK4) for N=10, 100, and 1000.}
    \label{fig: Figure 6}
\end{figure}

Naturally, we expect the maximum error to decrease with increasing $N$ until it reaches machine precision, and this is exactly what happens in Figure 8.
\begin{figure}[h!]
    \centering
    \includegraphics[width=1\linewidth]{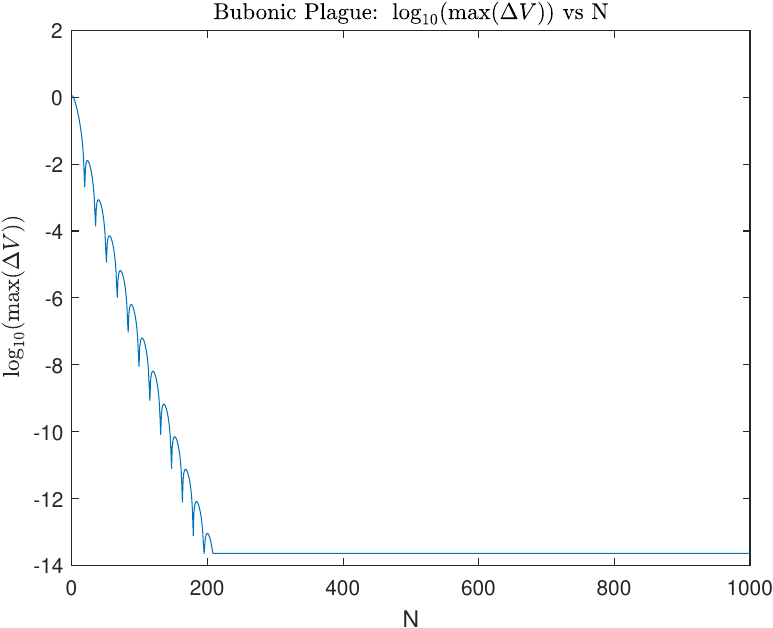}
    \caption{For the 1966 Bubonic Plague, this is a plot of the log of the maximum error of V (relative to RK4) versus the value of N.}
    \label{fig: Figure 7}
\end{figure}

Also, the values of the $V$ coefficients are strongly correlated with the calculated error in the case of the Bubonic plague as seen in Figure 9. This makes sense since the largest value $y$ can be is $1$, and so, the maximum next order correction is assumed to be roughly equal to the absolute value of the next coefficient. Note that there is also some further discussion in Appendix N on why the coefficients take on the generic convergent graphical behavior seen in Figures 9 and 12 as well as elsewhere in the results in Figures 17 and 19.
\begin{figure}[h!]
    \centering
    \includegraphics[width=1\linewidth]{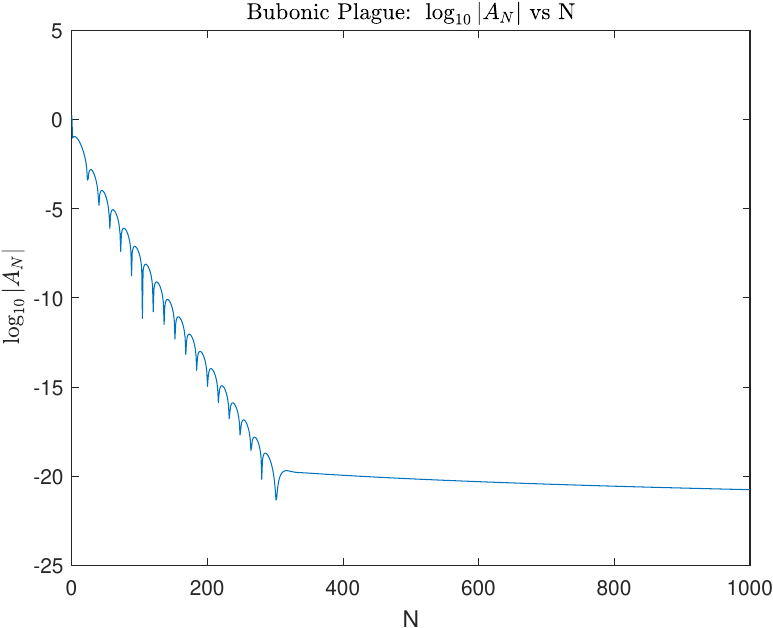}
    \caption{For the 1966 Bubonic Plague, this is a semilog plot of the absolute value of the coefficients $A_{N}$ (from the shifted exponential series of $V$) versus $N$.}
    \label{fig: Figure 8}
\end{figure}

In fact, the same previous comments about the 1966 Bubonic Plague can also be said for the Ebola scenario (Figures 10, 11, and 12), and in general, we will see a commonality on the behavior of scenarios with initial conditions located inside the HK region.
\begin{figure}[h!]
    \centering
    \includegraphics[width=1\linewidth]{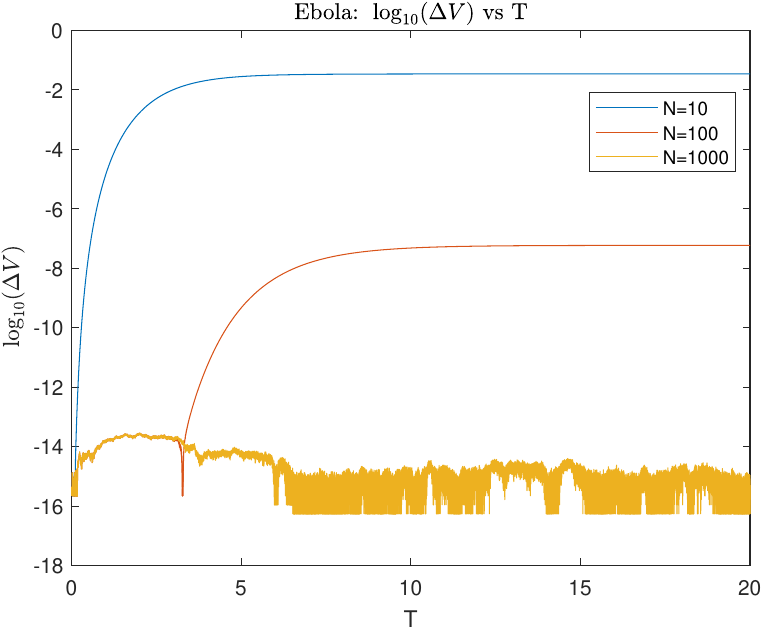}
    \caption{In this Ebola plot, we have the logs of the errors of V calculated (in comparison to RK4) for N=10, 100, and 1000.}
    \label{fig: Figure 9}
\end{figure}
\begin{figure}[h!]
    \centering
    \includegraphics[width=1\linewidth]{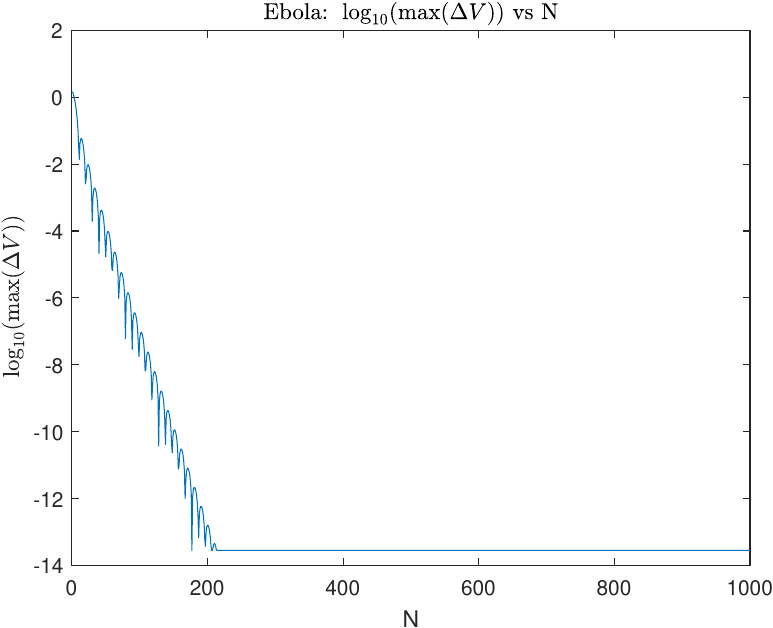}
    \caption{For Ebola, this is a plot of the log of the maximum error of V (relative to RK4) versus the value of N.}
    \label{fig: Figure 10}
\end{figure}
\begin{figure}[h!]
    \centering
    \includegraphics[width=1\linewidth]{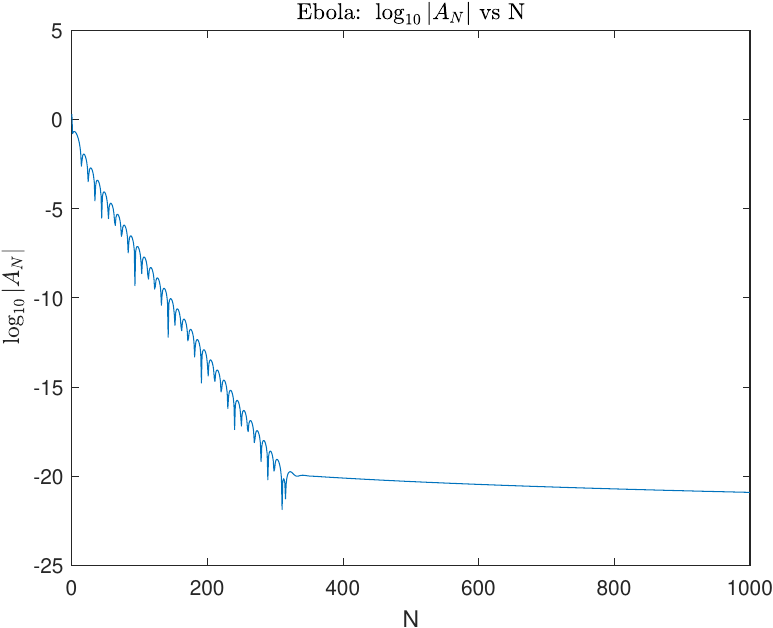}
    \caption{For Ebola, this is a semilog plot of the absolute value of the coefficients $A_{N}$ (from the shifted exponential series of $V$) versus N .}
    \label{fig: Figure 11}
\end{figure}

Next, there is the Covid-19 scenario in Japan which has a finite radius of convergence (in nondimensionalized time). Note that the divergence occurs shortly before all three graphs intersect for different chosen $N\boldsymbol{=} 10, 100, 1000$.
\begin{figure}[h!]
    \centering
    \includegraphics[width=1\linewidth]{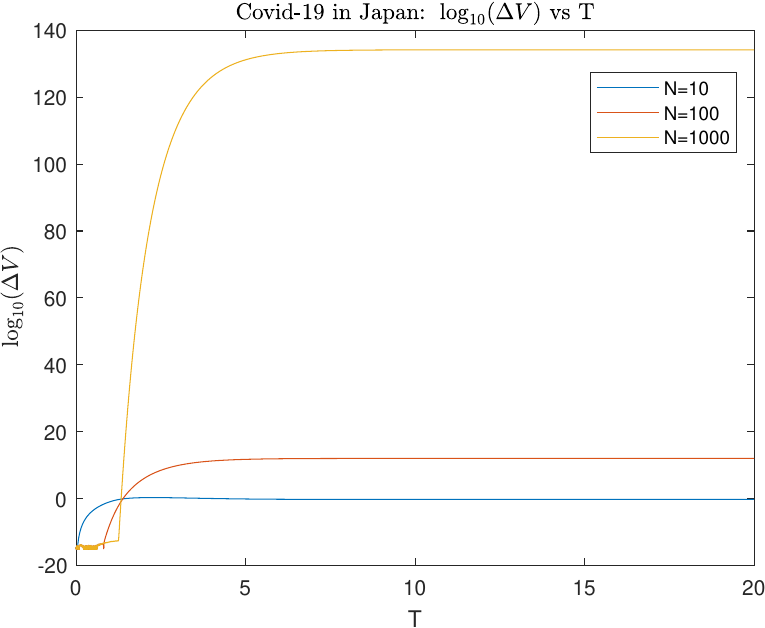}
    \caption{In this Covid-19 plot, we have the logs of the errors of V calculated (in comparison to RK4) for N=10, 100, and 1000.}
    \label{fig: Figure 12}
\end{figure}

\subsection{1966 Bubonic Plague: Mapping the Singularities}

Note that the singularities of a power series can be found at the locations of the zeroes of the reciprocal of the series in question. In particular, we would like to see what happens in the mapping of singularities from the direct time series to the shifted exponential series. Let us consider the case of the 1966 Bubonic Plague to illustrate this.

To find the series of the reciprocal, the following mathematical procedure can be done.
\begin{equation}
    \left(\sum_{n\boldsymbol{=} 0}^{\infty} b_{n}x^{n}\right)^{-1}\boldsymbol{=}\sum_{n\boldsymbol{=} 0} c_{n}x^{n},\,\,\, b_{0}\neq 0,\,\,\, c_{n>0}\boldsymbol{=}-\frac{1}{b_{0}}\sum_{j=1}^{n}b_{j}c_{n-j},\,\,\, c_{0}\boldsymbol{=}\frac{1}{b_{0}}
\end{equation}

Now, in the case of the direct time series for the Bubonic Plague, we would like to use this procedure to numerically find all the singularities that exist on the complex plane. Due to numerical sensitivity issues with MATLAB's roots function for relatively large degree polynomials, this was done with the variable-precision arithmetic function throughout the overall algorithm, and so, given how computationally expensive this makes it, there are practical limits with how large we can make $N$. With that in mind, the following singularities (in Table 1) which are closest to the origin were calculated for various values of $N$.
\begin{table}[h!]
    \centering
    \begin{tabular}{|c|c|c|}
    \hline
         $N$ & Singularities of $T$ Series& $\rho_{T}$\\
         \hline
        $10$ & $3.136416979332549\pm 4.716647383373361i$ & $5.664262882964409$\\
        \hline
        $18$ & $3.041750097587512\pm 4.457486116243105i$ & $5.396427163658719$\\
        \hline
        $32$ & $2.908760467592041\pm 4.316980448108888i$ & $5.205497828947840$\\
        \hline
        $56$ & $2.841660240593447\pm 4.233364593343863i$ & $5.098667345801946$\\
        \hline
        $100$ & $2.793295077275252\pm 4.181035559841080i$ & $5.028275623052675$\\
        \hline
        $178$ & $2.765280611466495\pm 4.150005924993069i$ & $4.986915483305295$\\
        \hline
        $214$ & $2.758820421556256\pm 4.142540491442002i$ & $4.977120838560420$\\
        \hline
        $316$ & $2.749120742536354\pm 4.131048556840228i$ & $4.962179665833893$\\
        \hline
    \end{tabular}
    \caption{This is a table of approximate singularities found from the zeroes of the polynomial of the reciprocal of the direct time series for the Bubonic plague example.}
    \label{tab: Table 1}
\end{table}

Next, doing the same procedure again with the shifted exponential series for the 1966 Bubonic Plague, we are able to obtain the following list of nearest singularity conjugate pairs (Table 2) for various values of $N$. From $N\boldsymbol{=} 178$ (or earlier), the real part of the nearest singularity pair appears to reasonably stabilize to around $1.048$, but the imaginary component continues to evolve towards a smaller number in the table. However, overall, the radius of convergence $\rho$ is stabilizing to around $1.04863$. In fact, according to the test, $|A_{1000}|^{-1/1000}\boldsymbol{=} 1.048926983566503$ which is in reasonable agreement (albeit not perfect) with our estimation from the singularities.

\begin{table}[h!]
    \centering
    \begin{tabular}{|c|c|c|}
    \hline
      $N$   & Singularities of Shifted Exponential Series& $\rho_{y}$ \\
      \hline
       $10$  & $1.056842793206091\pm 0.485710037224205i$ & $1.163112604098155$ \\
       \hline
       $18$  & $1.036143795072208\pm 0.283362835680417i$ & $1.074192003652739$\\
       \hline
       $32$  & $1.037765731959409\pm 0.175282863377030i$ & $1.052464629630326$ \\
       \hline
       $56$  & $1.043755632992874\pm 0.107370586044249i$ & $1.049263677133561$\\
       \hline
       $100$ & $1.046872197250276\pm 0.062602875686176i$ & $1.048742350351029$ \\
       \hline
       $178$ & $1.048034587685202\pm 0.035965767735149i$ & $1.048651530983133$ \\
       \hline
       $214$ & $1.048212100828620\pm 0.030060161212960i$ & $1.048643038224017$ \\
       \hline
       $316$ & $1.048434308295043\pm 0.020514745692475i$ & $1.048634995411145$ \\
       \hline
    \end{tabular}
    \caption{This is a table of approximate singularities found from the zeroes of the polynomial of the reciprocal of the $y$ series for the Bubonic plague example.}
    \label{tab: Table 2}
\end{table}

This brings us to the next detail to confirm. Does the closest singularity in the $T$ domain map to the closest singularity in the $y$ domain? Recall in the $T$ domain that we have $T_{b}\boldsymbol{=} \tau\pm i\omega\boldsymbol{=} 2.749120742536354\pm 4.131048556840228i$. A mapping from one singularity to another will occur, but surprisingly, the result of the mapping is not quite in complete agreement compared to what we obtained for the nearest singularity in the $y$ domain.
\begin{subequations}
    %\begin{equation}
    \begin{align}
    y_{b} &\boldsymbol{=}1-e^{\lambda T_{b}}\boldsymbol{=} (1-e^{\lambda\tau_{b}}\cos(\lambda\omega_{b}))\mp i e^{\lambda\tau_{b}}\sin(\lambda\omega_{b})\\
    & \boldsymbol{=} 1.122464903359770\pm 0.218081906713538i
    \end{align}
    %\end{equation}
    \begin{equation}
        \rho_{y,b}=|y_{b}|\boldsymbol{=}\sqrt{1+e^{2\lambda\tau_{b}}-2e^{\lambda\tau_{b}}\cos(\lambda\omega_{b})}\boldsymbol{=} 1.143454055618445
    \end{equation}
\end{subequations}

Given the tension that exists from the numerical values seen in both the root test and Table 2 versus equation (47b), it is worth noting that a comparison of the root test for $N\boldsymbol{=} 1000$ and $N\boldsymbol{=} 1001$ suggests an uncertainty in the radius of convergence $\Delta\rho_{bubonic}\boldsymbol{=} |A_{1000}^{-1/1000}-A_{1001}^{-1/1001}|\approx 4.794061345503309\times 10^{-5}$.

Unlike the case of the direct time series, it can be noted that the variable precision arithmetic does not have any noticeable effect on the calculated closest singularities of the $y$ series. So, to further clarify this issue regarding the tension between the root test and Table 2, another table of singularities (for the $y$ series) is calculated with larger values of $N$ without using the variable precision arithmetic. This has the advantage of allowing us to map from the $y$ series closest singularity to the $T$ series closest singularity for a much larger $N$ value.

\begin{table}[h!]
    \centering
    \begin{tabular}{|c|c|c|}
    \hline
      $N$   & Singularities of Shifted Exponential Series& $\rho_{y}$ \\
      \hline
       $400$  & $1.048507030362304\pm 0.016261862764783i$ & $1.048633129793141$ \\
       \hline
       $500$  & $1.048551119968439\pm 0.013042852480336i$ & $1.048632236385994$\\
       \hline
       $600$  & $1.048575315780594\pm 0.010887631559347i$ & $1.048631838819204$ \\
       \hline
       $700$  & $1.048590010475455\pm 0.009343675576619i$ & $1.048631639014480$\\
       \hline
       $800$ & $1.048599783127986\pm 0.008184349853659i$ & $1.048631722178280$ \\
       \hline
       $900$ & $1.048615701318654\pm 0.007362892510257i$ & $1.048641550406110$ \\
       \hline
       $1000$ & $1.037911784988708\pm 0.120142550117864i$ & $1.044842143946764$ \\
       \hline
       $1100$ & $1.035149423509403\pm 0.108874252051280i$ & $1.040859227634350$ \\
       \hline
       $1200$ & $1.032233792801933\pm 0.104922871742221i$ & $1.037552606867189$ \\
       \hline
       $1300$ & $1.029749762322831\pm 0.101594444326236i$ & $1.034749246978167$ \\
       \hline
       $1400$ & $1.028040214339160\pm 0.094141154453040i$ & $1.032341629142338$ \\
       \hline
       $1500$ & $1.026134314996485\pm 0.092001883451537i$ & $1.030250444781236$ \\
       \hline
       $1600$ & $1.026031112794913\pm 0.069980225083115i$ & $1.028414836690842$ \\
       \hline
       $1700$ & $1.023578688961388\pm 0.080932643408311i$ & $1.026773307631715$ \\
       \hline
       $1800$ & $1.021301957127695\pm 0.090615227338533i$ & $1.025314004126772$ \\
       \hline
       $1900$ & $1.020290836269810\pm 0.085701266977750i$ & $1.023883830196444$ \\
       \hline
       $2000$ & $1.011907849593462\pm 0.148133486304518i$ & $1.022693026197791$ \\
       \hline
    \end{tabular}
    \caption{This is a table of approximate singularities found from the zeroes of the polynomial of the reciprocal of the $y$ series for the Bubonic plague example with $N$ varying from $N\boldsymbol{=} 400$ to $N\boldsymbol{=} 2000$.}
    \label{tab: Table 3}
\end{table}

The last row gives a radius of convergence equal to $1.022693026197791$ which is in reasonable agreement with the root test $A_{2000}^{-1/2000}\boldsymbol{=} 1.024885997225322$. Note though the ratio test gives $A_{2000}/A_{2001}\boldsymbol{=} 1.001005121126731$. So, it is clear there is a slow convergence between the singularity method, root test, and ratio test.

\begin{figure}[h!]
    \centering
    \includegraphics[width=1\linewidth]{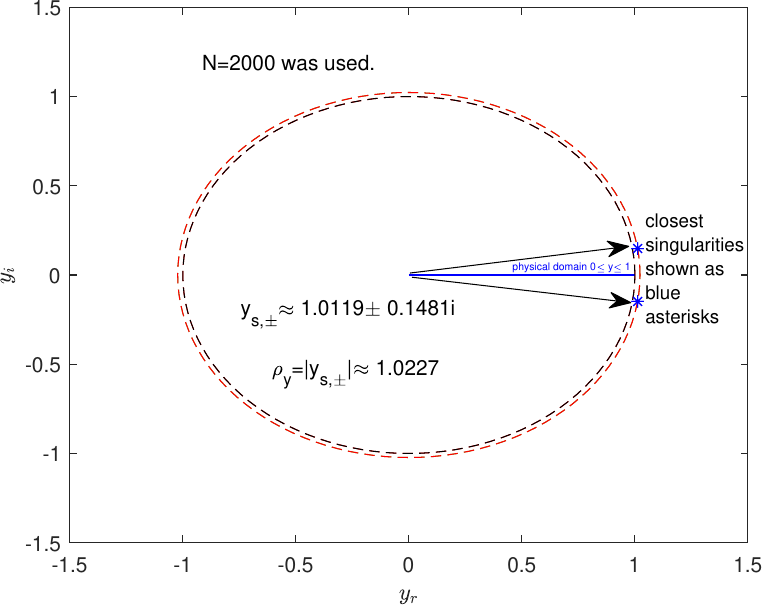}
    \caption{For the 1966 Bubonic Plague scenario, the closest singularities are illustrated and are shown to be outside of the physical domain in the $y$ domain with $N=2000$.}
    \label{fig:Y Radius from Singularities}
\end{figure}

Mapping the closest singularity (with $N\boldsymbol{=}2000$) from the $y$ domain to the $T$ domain results in the following. Furthermore, as we will see, this will again give a radius of convergence very close to $5$ in the $T$ domain which is in our expectations for the Bubonic plague based on what we saw in Figure 3. In doing this calculation, we are finally able to that we can get a set of numbers more reasonably in agreement with both the root tests and with what is seen in the mapping compared to the graphical expectation.
\begin{equation*}
    T_{c}\boldsymbol{=}\frac{1}{\lambda}\ln(1-y_{c})\approx\frac{1}{\lambda}\ln(-0.011907849593462\mp 0.148133486304518i)
\end{equation*}
\begin{equation*}
    \approx \frac{1}{\lambda}(-1.906420921805951\mp 1.651009786760101i)
\end{equation*}
\begin{equation}
    \approx 3.781820350115900\pm 3.275154158450535i
\end{equation}
\begin{equation}
    \rho_{T,c}\boldsymbol{=} |T_{c}|\approx 5.002879163258551
\end{equation}

\begin{figure}[h!]
    \centering
    \includegraphics[width=1\linewidth]{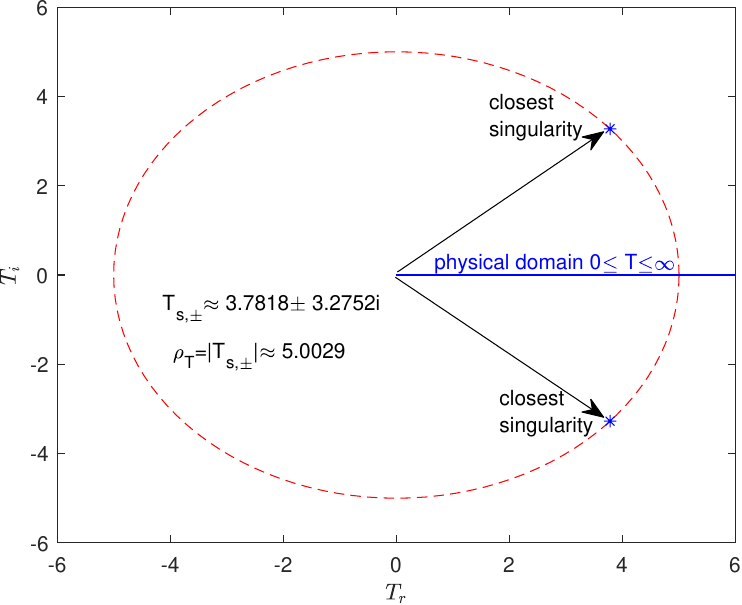}
    \caption{Through a mapping from the $y$ domain with $N=2000$, the singularities are found in the $T$ domain.}
    \label{fig: T radius from mapped Singularities}
\end{figure}

\subsection{Some Further Examples from the Radius Contour Plot}

Looking at some other examples from the Hershey Kiss region, the chosen case of $(\tilde{S_{0}}, \tilde{I_{0}})\boldsymbol{=} (2, 1)$ is an interesting case to observe in Figures 16 and 17. Take note of the period group pattern behavior in its maximum error and coefficient plots. Just as in the case of the Bubonic Plague and Ebola coefficient semilog plots, there is an eventual flattening out into a tail or more precisely a line with nearly zero slope.
\begin{figure}[h!]
    \centering
    \includegraphics[width=1\linewidth]{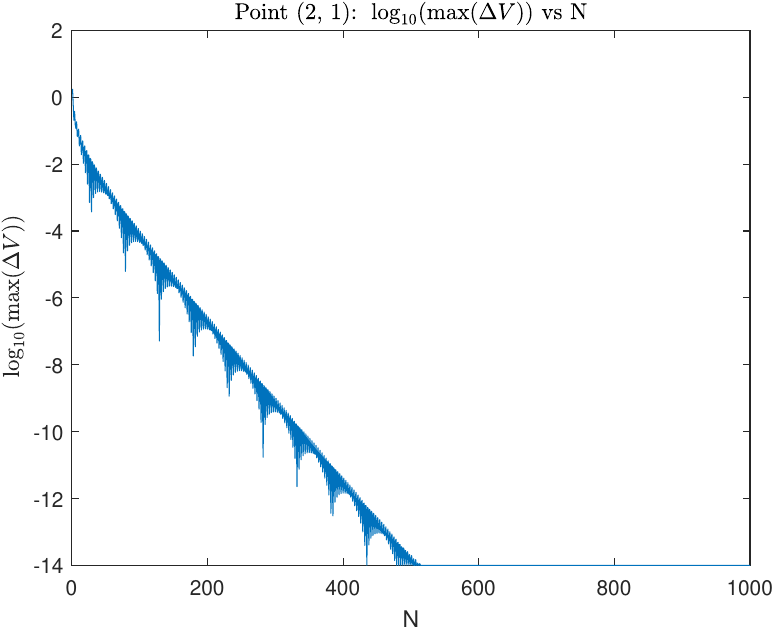}
    \caption{This is a semilog plot of the maximum error of the shifted exponential series when $(\tilde{S_{0}}, \tilde{I_{0}})\boldsymbol{=} (2, 1)$.}
    \label{fig: Figure 13}
\end{figure}

\begin{figure}[h!]
    \centering
    \includegraphics[width=1\linewidth]{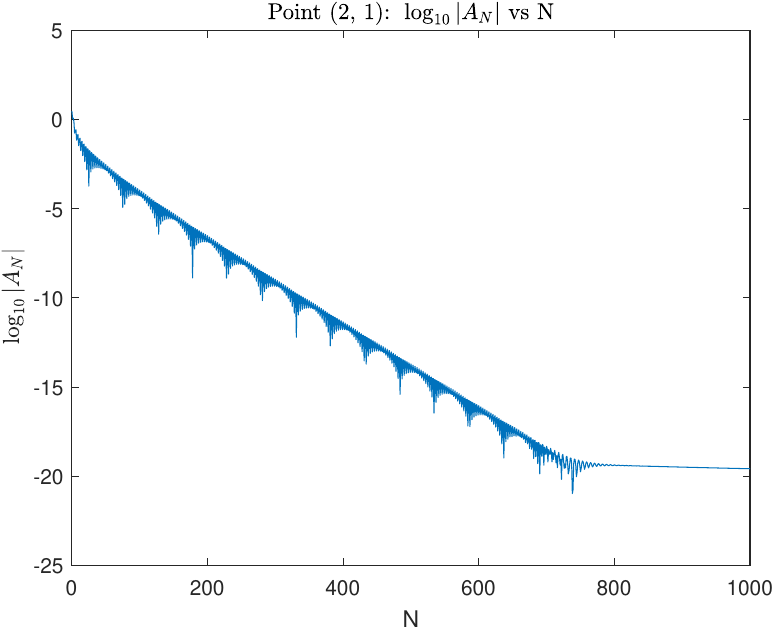}
    \caption{This is a semilog plot of the absolute values of the shifted exponential series coefficients when $(\tilde{S_{0}}, \tilde{I_{0}})\boldsymbol{=} (2, 1)$.}
    \label{fig: Figure 14}
\end{figure}

Compared to the Bubonic Plague and Ebola examples, note that this example does require considerably more calculations of coefficients before it reaches a set of steady-state coefficients. Next, let us try to select several examples near the likely maximum $\tilde{S_{0}}$ with a very small $\tilde{I_{0}}$.

If we look carefully at the example of $(\tilde{S_{0}}, \tilde{I_{0}})\boldsymbol{=} (2.47, 0.001)$ in Figure 18, we initially see behavior that appears to show full convergence during the first thousand coefficients. However, carrying out calculations to $100$ thousand coefficients, it instead becomes abundantly clear that the series here is divergent with some radius of convergence less than $1$. Of course, since the Hershey Kiss region is based on the $1000$th coefficient, this example shows up falsely as convergent. It is also clear that this example has an optimal truncation (for pretending full convergence) at the $3509$th coefficient. In general, this behavior is also consistent with the toy model in Appendix N.
\begin{figure}[h!]
    \centering
    \includegraphics[width=1\linewidth]{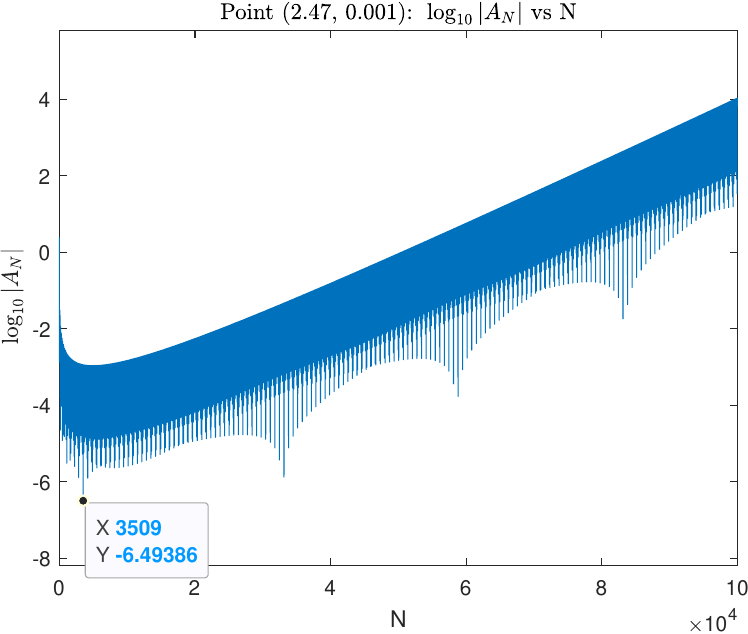}
    \caption{This is a semilog plot of the absolute values of the shifted exponential coefficients when $(\tilde{S_{0}}, \tilde{I_{0}})\boldsymbol{=} (2.47, 0.001)$.}
    \label{fig: Figure 15}
\end{figure}

Next, let us consider the example of $(\tilde{S_{0}}, \tilde{I_{0}})\boldsymbol{=} (2.45, 0.001)$ in Figure 19, and in fact, this case is clearly fully convergent although it takes roughly $240$ thousand coefficients for us to obtain a steady-state tail the semilog plot of the coefficients.
\begin{figure}[h!]
    \centering
    \includegraphics[width=1\linewidth]{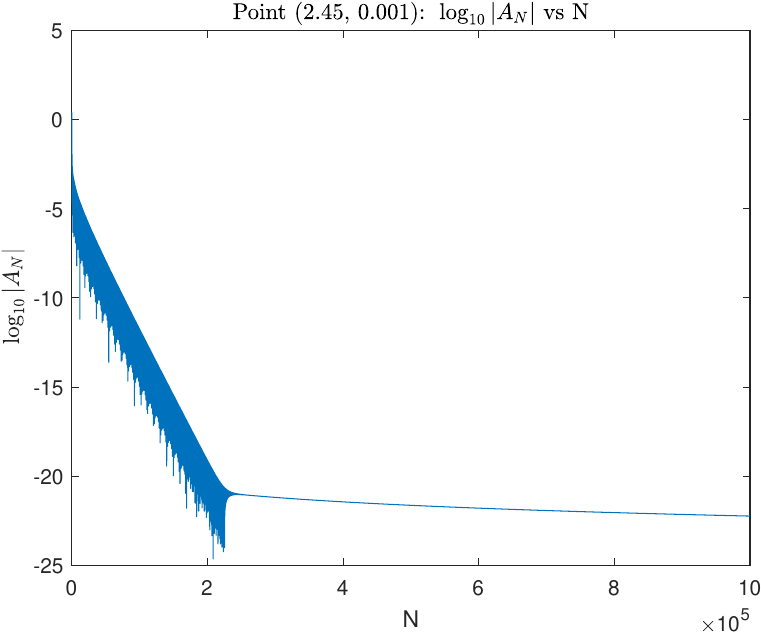}
    \caption{This is a semilog plot of the absolute values of the coefficients when $(\tilde{S_{0}}, \tilde{I_{0}})\boldsymbol{=} (2.45, 0.001)$.}
    \label{fig: Figure 16}
\end{figure}

\section{Discussion and Future Work}
\subsection{Why the Shifted Exponential Gauge Works}

The shifted exponential gauge has the major advantage of preserving the simplicity of the original initial condition thereby avoiding the computationally expensive adventure of solving for $A_{1}$ in the straight exponential gauge. Due to what was observed with the mapping of singularities related to the 1966 Bubonic Plague point, it can be seen that in general the nearest singularity pair in the $T$ domain does roughly map to the nearest singularity pair in the $y$ domain. In the author's opinion, in comparison to the shifted exponential series, the direct time series is probably more slowly convergent, and therefore, it likely requires many more terms in order for the singularity root finding to result in a more accurate mapping. This would also explain why the uncertainty in the Bubonic radius of convergence is not large enough to account for the tension between the different numbers although the later anaylysis with $N=2000$ largely resolved the tension through the reverse mapping.

\subsection{Future Work}

We want to use further re-summing techniques to increase the region of convergence for the shifted exponential gauge. Since it appears that all examples in the divergent region have various kinds of complicated repeating sign patterns that go on forever with the sequence of coefficients, a substitution of the following form where $\tilde{\rho}$ is a free parameter, has a good chance of doing so.

\begin{equation}
    z\boldsymbol{=}\frac{y}{y+\tilde{\rho}}
\end{equation}

As part of the goal of increasing the area of full physical convergence, we would like to better understand the nontrivial mapping of the singularities. This is likely responsible for the shape of the HK region, and note that it can also be observed, that as a function of $\tilde{I_{0}}$, the maximum value $\tilde{S}_{0, max}$ appears to be possibly a Gaussian function, or at the very least, it is very similar to a Gaussian curve. Note that this is the boundary of convergence by definition. This also includes a critical number $\tilde{S}_{0, c}$ effectively located at $(\tilde{S}_{0, c}, 0)$ in the off-axis limit $\tilde{I_{0}}\rightarrow 0$.

Although it is counterintuitive, if one considers how $\tilde{I}$ varies as a function of $\tilde{S}$, this result appears to be made possible by the on-axis trajectory behavior not being a limit of the off-axis trajectory. In short, when $\tilde{S_{0}}>1$ and $\tilde{I_{0}}=0$, we have an unstable equilibrium. To explain this point in more detail, as soon as $\tilde{I_{0}}$ is greater than zero even when it is extremely small ($\tilde{I_{0}}<<1$), $\tilde{I}$ then evolves towards a maximum which happens when $\tilde{S}\boldsymbol{=} 1$ occurs giving the maximum $\tilde{I}$ to be a function of $\tilde{S_{0}}$ in the limit of $\tilde{I_{0}}\rightarrow 0$. This maximum $\tilde{I}(T)$ can be arbitrarily large for large $\tilde{S_{0}}$ in this limit, but since this peak suddenly disappears when actually choosing $\tilde{I_{0}}\boldsymbol{=} 0$ instead of having $\tilde{I_{0}}\rightarrow 0$, the existence of $(\tilde{S}_{0, c}, 0)$ in the corner of the HK region begins to seem much more plausible. So, clarifying this issue through this line of thought should occur in any future work on the shape of the fully convergent parameter space.

This is a strange result since one might assume that there should be an infinite radius of curvature at that point. However, some initial asymptotic analysis in the appendix has shown why the radius of curvature should actually be $1$. So, the next question then is if further asymptotic analysis can shed further light on this critical number and the shape of the boundary of convergence.

\section{Conclusions}
In conclusion, a convergent power series solution is obtained for the SIR model, using an asymptotically motivated gauge function. For certain choices of model parameter values, the series converges over the full physical domain (i.e., for all positive time), and for other choices, this fails to occur. Consequently, perhaps, due to the more complicated dynamics of singularity mapping, there is a boundary of convergence within the parameter space, and in general, the radius of convergence is clearly a function of the initial nondimensionalized susceptible and infected populations. Further resummation techniques can be applied to increase the area of the convergent parameter space.

\section{Appendices}
\subsection{Appendix A: Simplifying the Problem}
If we define the following, a simplification results from equation (\ref{SWB}).
\begin{equation}
    \xi\boldsymbol{=} S e^{-\beta/\alpha}
    \label{xiEQN}
\end{equation}

Essentially, the square has been completed with the above substitution.
\begin{equation}
    \frac{d\xi}{dt}\boldsymbol{=} r e^{\beta/\alpha}\xi^{2}-\alpha\xi\ln(\xi)
    \label{xi2}
\end{equation}

Next, we would like to choose a non-dimensional time.
\begin{equation}
    T\boldsymbol{=}\alpha t
\end{equation}

With this choice, equation (\ref{xi2}) now transforms into the following.
\begin{equation}
    \frac{d\xi}{dT}\boldsymbol{=} L\xi^{2}-\xi\ln(\xi)
    \label{xi3}
\end{equation}

Note that the combined parameter $L$ has the following definition.
\begin{equation}
    L\boldsymbol{=} \frac{r}{\alpha} e^{\beta/\alpha}
    \label{L}
\end{equation}

\subsection{Appendix B: Rewriting the Original Parameters}
The multiplying factor $e^{-\beta/\alpha}$ from equation (\ref{xiEQN}) can be rewritten as the following.
\begin{equation}
    e^{-\beta/\alpha}\boldsymbol{=} e^{-\ln(S_{0})+r(S_{0}+I_{0})/\alpha}\boldsymbol{=}\frac{1}{S_{0}}e^{r(S_{0}+I_{0})/\alpha}
\end{equation}

The combined parameter $L$ can also be rewritten as the following.
\begin{equation}
    L\boldsymbol{=}\frac{r S_{0}}{\alpha}e^{-r(S_{0}+I_{0})/\alpha}
    \label{L2}
\end{equation}

In particular, the form of equation (\ref{L2}) suggests it will be convenient to define the following two parameters.
\begin{equation}
    \tilde{S_{0}}\boldsymbol{=}\frac{r S_{0}}{\alpha}
\end{equation}
\begin{equation}
    \tilde{I_{0}}\boldsymbol{=}\frac{r I_{0}}{\alpha}
\end{equation}

More generally, it might also be convenient to define the following as well.
\begin{equation}
    \tilde{S}\boldsymbol{=}\frac{r S}{\alpha}
\end{equation}
\begin{equation}
    \tilde{I}\boldsymbol{=}\frac{r I}{\alpha}
\end{equation}

Consequently, the combined parameter can be rewritten as the following.
\begin{equation}
    L\boldsymbol{=}\tilde{S_{0}}e^{-\tilde{S_{0}}-\tilde{I_{0}}}
    \label{L3}
\end{equation}

Also, the multiplying factor $e^{-\beta/\alpha}$ can be rewritten as the following.
\begin{equation}
    e^{-\beta/\alpha}\boldsymbol{=}\frac{1}{S_{0}}e^{\tilde{S_{0}}+\tilde{I_{0}}}
\end{equation}

As a result, equation (\ref{xiEQN}) tells us that the initial value of $\xi$ will be the following.
\begin{equation}
    \xi_{0}\boldsymbol{=} e^{\tilde{S_{0}}+\tilde{I_{0}}}
\end{equation}

Consequently, the combined parameter can also be written as the following.
\begin{equation}
    L\boldsymbol{=}\frac{\tilde{S_{0}}}{\xi_{0}}
    \label{L4}
\end{equation}

As an additional aside, we can obtain the following from equation (\ref{xiEQN}) with our new definitions.
\begin{equation}
    \tilde{S}\boldsymbol{=}\left(\frac{r}{\alpha}\right)\xi e^{\beta/\alpha}\boldsymbol{=}\tilde{S_{0}}\xi e^{-\tilde{S_{0}}-\tilde{I_{0}}}\boldsymbol{=} L\xi
    \label{Stilde2}
\end{equation}

Finally, combining equations (\ref{L4}) and (\ref{Stilde2}) gives the following.

\begin{equation}
    L\boldsymbol{=}\frac{\tilde{S}}{\xi}\boldsymbol{=}\frac{\tilde{S_{0}}}{\xi_{0}}
\end{equation}

Note that one could have also made an intermediate equation before introducing equation (\ref{xi3}) in which the dynamical function was instead $\tilde{S}$. This would give us the corresponding equation below.

\begin{equation}
    \frac{\tilde{S}}{dT}\boldsymbol{=}-(\tilde{S_{0}}+\tilde{I_{0}})\tilde{S}+\tilde{S}^{2}-\tilde{S}\ln(\tilde{S})
\end{equation}

\subsection{Appendix C: Enter the Lambert Function}
Next, we wish to solve for the steady state value value of equation (\ref{xi3}). Naturally, this will be $\xi_{\infty}$ at $T\boldsymbol{=}\infty$ which leads to $d\xi/dT\boldsymbol{=} 0$. Solving for $\xi_{\infty}$, we have the following.
\begin{equation}
    L\xi_{\infty}^{2}-\xi_{\infty}\ln(\xi_{\infty})\boldsymbol{=} 0
\end{equation}

$\xi_{\infty}\boldsymbol{=} 0$ is the trivial solution which can be discarded.
\begin{equation}
    L\xi_{\infty}-\ln(\xi_{\infty})\boldsymbol{=} 0
    \label{xiInfinityEQN}
\end{equation}

Manipulating the equation further, we arrive at the following.
\begin{equation}
    -L\xi_{\infty}e^{-L\xi_{\infty}}\boldsymbol{=} -L
    \label{xiL}
\end{equation}

Equation (\ref{xiL}) can be solved with the Lambert functions.
\begin{equation}
    -L\xi_{\infty}\boldsymbol{=} W_{0}(-L),\, W_{-1}(-L)
\end{equation}
\begin{equation}
    \xi_{\infty}\boldsymbol{=}-\frac{1}{L}W_{0}(-L),\, -\frac{1}{L}W_{-1}(-L)
\end{equation}

Note that both Lambert derived solutions satisfy the following identity.
\begin{equation}
    \ln(\xi e^{-L\xi})\boldsymbol{=}0
\end{equation}

Ultimately, we must narrow this down to only one solution, and the criteria for deciding this will be to require that our steady state solution is stable. Before we do this, let us take note of several properties for further clarification of this discussion with an initial focus on the maximum value of the combined parameter (which occurs at $\tilde{S_{0}}=1$).
\begin{equation}
    L_{max}\boldsymbol{=} max(\tilde{S_{0}}e^{-\tilde{S_{0}}})e^{-\tilde{I_{0}}}\boldsymbol{=} e^{-1-\tilde{I_{0}}}
\end{equation}

As a result, the following inequality is also true.
\begin{equation}
    0\leq L\leq\frac{1}{e}
\end{equation}

Of course, note that the absolute minimum case $L\boldsymbol{=} 0$ occurs whenever $\tilde{S_{0}}\boldsymbol{=} 0$, and the absolute maximum $L\boldsymbol{=} 1/e$ occurs when $\tilde{S_{0}}\boldsymbol{=} 1$ and $\tilde{I_{0}}\boldsymbol{=} 0$. Some further discussion of these special cases will be made in a later section.

Next, let us take note of an equality of the Lambert functions when $-1/e\leq x\leq 0$.
\begin{equation}
    W_{-1}(x)\leq W_{0}(x)\leq 0
    \label{WEQ}
\end{equation}

Also, for the same domain input of $x$, the following two inequalities are true as well.
\begin{equation}
    -1\leq W_{0}(x)\leq 0
\end{equation}
\begin{equation}
    W_{-1}(x)\leq -1
\end{equation}

By multiplication of $-1/L$ onto equation (\ref{WEQ}), then the following is trivially true with $x\boldsymbol{=}-L$.
\begin{equation}
    0\leq -\frac{1}{L}W_{0}(-L)\leq -\frac{1}{L}W_{-1}(-L)
\end{equation}

Therefore, to determine the question of stability, we must consider sign of the derivative of $d\xi/dT$ with respect to $\xi$ inside the different intervals of
    $0\leq\xi\leq -\frac{1}{L}W_{0}(-L)$, $-\frac{1}{L}W_{0}(-L)\leq\xi\leq -\frac{1}{L}W_{-1}(-L)$, and $-\frac{1}{L}W_{-1}(-L)\leq\xi$.

\begin{equation}
    \frac{d}{d\xi}\left[\frac{d\xi}{dT}\right]\boldsymbol{=}2L\xi-\ln(\xi)-1\boldsymbol{=} L\xi-\ln(\xi e^{-L\xi})-1
\end{equation}

Choosing a value $\xi\boldsymbol{=}\epsilon$ which is arbitrarily close to zero, it easy to show the derivative in the first interval must be positive.
\begin{equation}
    \lim_{\xi \to 0^{+}} \frac{d}{d\xi}\left[\frac{d\xi}{dT}\right]\boldsymbol{=}\lim_{\xi \to 0^{+}} \left[2L\xi+\ln(1/\xi)-1\right]\boldsymbol{=}+\infty>0
\end{equation}

So, any arbitrarily small positive $\xi$ will move away from the zero solution showing that $\xi\boldsymbol{=} 0$ is an unstable zero. Next, let us consider the zero provided through the Lambert function $W_{0}(x)$.
\begin{equation}
    \frac{d}{d\xi}\left[\frac{d\xi}{dT}\right]_{\xi\boldsymbol{=}-\frac{1}{L}W_{0}(-L)}\boldsymbol{=}\left[L\xi-\ln(\xi e^{L\xi})-1\right]_{\xi\boldsymbol{=}-\frac{1}{L}W_{0}(-L)}\boldsymbol{=}-W_{0}(-L)-1\leq 0
\end{equation}

Therefore, $\xi\boldsymbol{=}-\frac{1}{L}W_{0}(-L)$ is a stable zero. Finally, let consider the zero provided through the other Lambert function $W_{-1}(x)$.

\begin{equation}
    \frac{d}{d\xi}\left[\frac{d\xi}{dT}\right]_{\xi\boldsymbol{=}-\frac{1}{L}W_{-1}(-L)}\boldsymbol{=}\left[L\xi-\ln(\xi e^{-L\xi})-1\right]_{\xi\boldsymbol{=}-\frac{1}{L}W_{-1}(-L)}\boldsymbol{=}-W_{-1}(-L)-1\geq 0
\end{equation}

So, the final case of $\xi\boldsymbol{=}-\frac{1}{L}W_{-1}(-L)$ is an unstable zero. Therefore, we are left with only one choice for a stable zero.

\begin{equation}
    \xi_{\infty}\boldsymbol{=}-\frac{1}{L}W_{0}(-L)
\end{equation}

\subsection{Appendix D: Validity of Initial Values}
Since the expected result is assumed to be that $S(t)$ is a decreasing monotonic function, it should also be the case for for $\xi(T)$ as well. However, the previous analysis tells us this will only occur when $-\frac{1}{L}W_{0}(-L)\leq \xi_{0}\leq -\frac{1}{L}W_{-1}(-L)$.

Fortunately, this potential contradiction can be resolved because the combined parameter $L$ is a function of the initial values $\tilde{S_{0}}$ and $\tilde{I_{0}}$, but we can make the argument by simply by having $L$ be a function of $\tilde{S_{0}}$ and $\xi_{0}$. In other words, $L=L(\tilde{S_{0}}, \tilde{I_{0}})$ and $L=L(\tilde{S_{0}}, \xi_{0})$ are both true.

\begin{equation}
    \left[\frac{d\xi}{dT}\right]_{T\boldsymbol{=} 0}\boldsymbol{=} L(\tilde{S_{0}}, \xi_{0})\xi_{0}^{2}-\xi_{0}\ln(\xi_{0})\boldsymbol{=}\left(\frac{\tilde{S_{0}}}{\xi_{0}}\right)\xi_{0}^{2}-\xi_{0}\ln\left(e^{\tilde{S_{0}}+\tilde{I_{0}}}\right)\boldsymbol{=}\tilde{S_{0}}\xi_{0}-\xi_{0}\left(\tilde{S_{0}}+\tilde{I_{0}}\right)\boldsymbol{=}-\tilde{I_{0}}\xi_{0}\leq 0
    \label{initialRate}
\end{equation}

Consequently, by simple substitution of appropriate identities, we can see that the derivative is always negative. Therefore, the function $\xi$ is always a decreasing monotonic function as required. Consequently, this also tells us that $-\frac{1}{L}W_{0}(-L)\leq \xi_{0}\leq -\frac{1}{L}W_{-1}(-L)$ is always true.

\subsection{Appendix E: A Couple of Special Cases}
Before solving the problem in general, let us consider a couple of special cases. First, let us consider the case of $\tilde{I_{0}}\boldsymbol{=} 0$. Since this special case leads to $\tilde{S_{0}}\boldsymbol{=}\ln(\xi_{0})$, this gives us the differential equation below.

\begin{equation}
    \frac{d\xi}{dT}\boldsymbol{=}\left(\frac{\ln(\xi_{0})}{\xi_{0}}\right)\xi^{2}-\xi\ln(\xi)
    \label{xi4}
\end{equation}

Simply inserting a constant $\xi\boldsymbol{=}\xi_{0}$ into equation (\ref{xi4}) automatically shows this is a general solution. We could have also shown this to be the only case satisfying the initial condition by substituting $\tilde{I_{0}}\boldsymbol{=} 0$ into equation (\ref{initialRate}). Also, note that in the more specific case of $L\boldsymbol{=} 1/e$, then implies $\tilde{I_{0}}\boldsymbol{=} 0$ and $\tilde{S_{0}}\boldsymbol{=} 1$ thereby giving us the constant solution $\xi\boldsymbol{=} e$ and corresponding $\tilde{S}\boldsymbol{=} 1$.

Another special case to consider is when $\tilde{S_{0}}\boldsymbol{=} 0$. Note that this results in the combined parameter becoming very simply $L\boldsymbol{=} 0$. Consequently, we have differential equation below.

\begin{equation}
    \frac{d\xi}{dT}\boldsymbol{=}-\xi\ln(\xi)
    \label{xi5}
\end{equation}

Solving equation (\ref{xi5}), it is simpler to use the substitution $\xi\boldsymbol{=} e^{V}$ which after cancellation of a factor gives us the following equation.

\begin{equation}
    \frac{dV}{dT}\boldsymbol{=}-V
\end{equation}

This has the following solution. Also, take note that $\ln(\xi_{0})\boldsymbol{=}\tilde{I_{0}}$ since $\tilde{S_{0}}\boldsymbol{=} 0$.

\begin{equation}
    V\boldsymbol{=} V_{0} e^{-T}\boldsymbol{=} \ln(\xi_{0}) e^{-T}\boldsymbol{=}\tilde{I_{0}} e^{-T}
    \label{44}
\end{equation}

\begin{equation}
    \xi\boldsymbol{=}\exp(\tilde{I_{0}} \exp(-T))
\end{equation}

However, since $\tilde{S_{0}}\boldsymbol{=} 0$, the combined parameter is $L\boldsymbol{=} 0$. As a result, we have the following by direct substitution in equation (\ref{Stilde2}).

\begin{equation}
    \tilde{S}\boldsymbol{=} L\xi\boldsymbol{=} 0\cdot \exp(\tilde{I_{0}} \exp(-T))\boldsymbol{=} 0
\end{equation}

\subsection{Appendix F: The First Few Coefficients of the Power Series of T}
Given the initial values, $A_{0}\boldsymbol{=}\tilde{S_{0}}+\tilde{I_{0}}$, $B_{0}\boldsymbol{=}\xi_{0}$, and $C_{0}\boldsymbol{=}\tilde{S_{0}}$, we can begin calculating the first few coefficients which gives us the following.

\begin{equation}
    A_{1}\boldsymbol{=}-\tilde{I_{0}}
\end{equation}

\begin{equation}
    B_{1}\boldsymbol{=}-\tilde{I_{0}}\xi_{0}
\end{equation}

\begin{equation}
    C_{1}\boldsymbol{=}-\tilde{I_{0}}\tilde{S_{0}}
\end{equation}

\begin{equation}
    A_{2}\boldsymbol{=}\frac{1}{2}\tilde{I_{0}}\left(1-\tilde{S_{0}}\right)
\end{equation}

\begin{equation}
    B_{2}\boldsymbol{=}\frac{1}{4}\tilde{I_{0}}\left(2\tilde{I_{0}}+1-\tilde{S_{0}}\right)\xi_{0}
\end{equation}

\begin{equation}
    C_{2}\boldsymbol{=}\frac{1}{4}\tilde{I_{0}}\left(2\tilde{I_{0}}+1-\tilde{S_{0}}\right)\tilde{S_{0}}
\end{equation}

\begin{equation}
    A_{3}\boldsymbol{=}\frac{1}{12}\tilde{I_{0}}\left(-2+3\tilde{S_{0}}+2\tilde{I_{0}}\tilde{S_{0}}-\tilde{S_{0}}^{2}\right)
\end{equation}

\begin{equation}
    B_{3}\boldsymbol{=}\frac{1}{12}\tilde{I_{0}}\left(-2\tilde{I_{0}}^{2}-5\tilde{I_{0}}+3\tilde{I_{0}}\tilde{S_{0}}+7\tilde{S_{0}}-2-\tilde{S_{0}}^{2}\right)\xi_{0}
\end{equation}

\begin{equation}
    C_{3}\boldsymbol{=}\frac{1}{12}\tilde{I_{0}}\left(-2\tilde{I_{0}}^{2}-5\tilde{I_{0}}+3\tilde{I_{0}}\tilde{S_{0}}+7\tilde{S_{0}}-2-\tilde{S_{0}}^{2}\right)\tilde{S_{0}}
\end{equation}

\begin{equation}
    A_{4}\boldsymbol{=}\frac{1}{48}\tilde{I_{0}}\left(2-5\tilde{S_{0}}-7\tilde{I_{0}}\tilde{S_{0}}+8\tilde{S_{0}}^{2}-2\tilde{I_{0}}^{2} \tilde{S_{0}}+3\tilde{I_{0}}\tilde{S_{0}}^{2}-\tilde{S_{0}}^{3}\right)
\end{equation}

\subsection{Appendix G: Dominant Balance of $\xi$}
Here, it will be noted that the dominant balance analysis of (\ref{xi3}) results in the same expected gauge function that we previously got in our analysis of (\ref{VEQN}). Once again, near the limit of $T\rightarrow\infty$, we will assume the function obeying equation (\ref{xi3}), is a combination of $\xi_{\infty}$ and a function $G$ which is much smaller than $\xi_{\infty}$.
\begin{equation}
    \xi\approx\xi_{\infty}+G
\end{equation}

Substituting this into equation (\ref{xi3}) now gives us the following.
\begin{equation}
    \frac{dG}{dT}\boldsymbol{=} L(\xi_{\infty}+G)^{2}-(\xi_{\infty}+G)\ln(\xi_{\infty}+G)
    \label{g1}
\end{equation}

Note that $\ln(1+x)\boldsymbol{=} x-x^{2}/2+x^{3}/3-...$, and so, the logarithm contribution gives the following.
\begin{equation*}
    (\xi_{\infty}+G)\ln(\xi_{\infty}+G)\boldsymbol{=}\xi_{\infty}(1+G/\xi_{\infty})(\ln(\xi_{\infty})+\ln(1+G/\xi_{\infty}))
\end{equation*}
\begin{equation*}
    \boldsymbol{=}\xi_{\infty}(1+G/\xi_{\infty})(\ln(\xi_{\infty})+G/\xi_{\infty}-(G/\xi_{\infty})^{2}/2+(G/\xi_{\infty})^{3}/3-...)
\end{equation*}

Keeping only first order terms, we get the following.
\begin{equation}
    (\xi_{\infty}+G)\ln(\xi_{\infty}+G)\boldsymbol{=} (\xi_{\infty}+G)\ln(\xi_{\infty})+G
\end{equation}
\begin{equation}
(\xi_{\infty}+G)^{2}\boldsymbol{=}\xi_{\infty}^{2}+2\xi_{\infty}G
\end{equation}

So, by simple substitution into equation (\ref{g1}), we obtain the following.
\begin{equation}
    \frac{dG}{dT}\boldsymbol{=} L\xi_{\infty}^{2}+2L\xi_{\infty}G-(\xi_{\infty}+G)\ln(\xi_{\infty})-G\boldsymbol{=} L\xi_{\infty}^{2}-\xi_{\infty}\ln(\xi_{\infty})+(2L\xi_{\infty}-\ln(\xi_{\infty})-1)G
    \label{g2}
\end{equation}

Note that the definition of $\xi_{\infty}$ as an asymptotic solution of $\xi$ satisfies both of the simple algebraic equations below.
\begin{equation}
    L\xi_{\infty}^{2}-\xi_{\infty}\ln(\xi_{\infty})\boldsymbol{=} 0
\end{equation}
\begin{equation}
     L\xi_{\infty}-\ln(\xi_{\infty})\boldsymbol{=} 0
\end{equation}

Consequently, equation (\ref{g2}) simplifies to the following result.
\begin{equation}
    \frac{dG}{dT}\boldsymbol{=} (L\xi_{\infty}-1)G
\end{equation}

Therefore, the simplest possible nontrivial solution in the dominant balance is the following.
\begin{equation}
    G\boldsymbol{=} c_{1} g\,\,\,\,\,\,\, g=e^{\lambda T}
    \label{g3}
\end{equation}
where
\begin{equation}
    \lambda\boldsymbol{=} L\xi_{\infty}-1\boldsymbol{=} -W_{0}(-L)-1
\end{equation}

\subsection{Appendix H: Some Analysis of the Straight Exponential Series}

With the definition of the $E_{n}$ and $F_{n}$ coefficients, (\ref{AgFinal}) transforms into the following.

\begin{equation}
    a_{1}+E_{2}a_{1}^2+E_{3}a_{1}^3+...\boldsymbol{=}\tilde{S_{0}}+\tilde{I_{0}}-\ln(\xi_{\infty})\boldsymbol{=}\ln(\tilde{S_{0}}/\tilde{S_{\infty}})
    \label{EgFinal}
\end{equation}

\begin{equation}
    F_{0}\boldsymbol{=}\tilde{S}_{\infty}
\end{equation}

\begin{equation}
    E_{1}\boldsymbol{=} 1
\end{equation}

\begin{equation}
    E_{n+1}\boldsymbol{=} \frac{1}{n(n+1)\lambda}\sum_{j\boldsymbol{=} 0}^{n-1} (j+1)E_{j+1}F_{n-j}
\end{equation}

\begin{equation}
    F_{n+1}\boldsymbol{=} ((n+1)\lambda+1)E_{n+1}
\end{equation}

\begin{equation}
    F_{1}\boldsymbol{=}\tilde{S}_{\infty}
\end{equation}

\begin{equation}
    E_{2}\boldsymbol{=}\frac{\tilde{S}_{\infty}}{2\lambda}
\end{equation}

\begin{equation}
    F_{2}\boldsymbol{=}\left(1+\frac{1}{2\lambda}\right)\tilde{S}_{\infty}
\end{equation}

\begin{equation}
    E_{3}\boldsymbol{=}\left(\frac{1}{6\lambda}+\frac{1}{12\lambda^{2}}\right)L\tilde{S}_{\infty}+\frac{\tilde{S}_{\infty}^{2}}{6\lambda^{2}}\boldsymbol{=}\left(\frac{1}{3}+\frac{1}{4\lambda}\right)\frac{\tilde{S}_{\infty}}{\lambda}
\end{equation}

\begin{equation}
    F_{3}\boldsymbol{=}\left(\frac{1}{3}+\frac{1}{4\lambda}\right)\left(3+\frac{1}{\lambda}\right)\tilde{S}_{\infty}\boldsymbol{=}\left(1+\frac{13}{12\lambda}+\frac{1}{4\lambda^{2}}\right)\tilde{S}_{\infty}
\end{equation}

The issue of convergence of this $g$ series is unclear within the parameter space of $(\tilde{S_{0}},\, \tilde{I_{0}})$, and furthermore, the issue of solving equation (\ref{EgFinal}) approximately for $A_{1}$ through root finding is quite computationally expensive for large $n$ while also raising complications of whether the correct root was selected.

\subsection{Appendix I: The First Few Coefficients of the Shifted Exponential Series}

Calculating the first few coefficients of the shifted exponential series gives us the following.

\begin{equation}
    A_{1}\boldsymbol{=}\frac{\tilde{I_{0}}}{\lambda}
\end{equation}

\begin{equation}
    C_{1}\boldsymbol{=}\frac{\tilde{I_{0}}\tilde{S_{0}}}{\lambda}
\end{equation}

\begin{equation}
    A_{2}\boldsymbol{=}\frac{\tilde{I_{0}}}{2\lambda}\left(1+\frac{1}{\lambda}-\frac{\tilde{S_{0}}}{\lambda}\right)
\end{equation}

\begin{equation}
    C_{2}\boldsymbol{=}\frac{\tilde{I_{0}}\tilde{S_{0}}}{2\lambda}\left(1+\frac{1}{\lambda}-\frac{\tilde{S_{0}}}{\lambda}+\frac{\tilde{I_{0}}}{\lambda}\right)
\end{equation}

\begin{equation}
    A_{3}\boldsymbol{=}\frac{\tilde{I_{0}}}{6\lambda}\left(2+\frac{3}{\lambda}-\frac{3\tilde{S_{0}}}{\lambda}+\frac{1}{\lambda^{2}}-\frac{2\tilde{S_{0}}}{\lambda^{2}}+\frac{\tilde{S_{0}}^{2}}{\lambda^{2}}-\frac{\tilde{I_{0}}\tilde{S_{0}}}{\lambda^{2}}\right)
\end{equation}

\subsection{Appendix J: An Asymptotic Analysis of Small $\tilde{I_{0}}$}

The lower right corner of the Hershey Kiss region can appear to be a counterintuitive result when one considers that being exactly on the $\tilde{S_{0}}$ axis, one should expect a constant quantity and therefore an infinite radius of convergence. However, the off-axis behavior even for arbitrarily very small $\tilde{I_{0}}$ indicates a finite radius leading to two obvious questions. Why does this happen? Why does there appear to be a critical value of $\tilde{S}_{0, max}\approx (2.47+2.45)/2=2.46$ in the limit of $\tilde{I_{0}}\rightarrow 0$ where the convergent parameter space comes to an end?

To answer the first question, the approach of asymptotics will be taken with an expansion of powers of small $\tilde{I_{0}}<<1$. Consider the following expansion.

\begin{equation}
    V\boldsymbol{=}\tilde{S_{0}}+\tilde{I_{0}}H_{1}(y)+\tilde{I_{0}}^{2}H_{2}(y)+\tilde{I_{0}}^{3}H_{3}(y)+...
\end{equation}

Substituting this expansion into equation (\ref{VEQNy}) while recalling the definition of $L$ in terms of $\tilde{S_{0}}$ and $\tilde{I_{0}}$ before finally focusing on the first order terms, the following equation of $H_{1}$ can then be obtained.

\begin{equation}
    \frac{dH_{1}}{dy}\boldsymbol{=}\frac{(1-\tilde{S_{0}})H_{1}+\tilde{S_{0}}}{\tilde{\lambda}(1-y)}
    \label{127}
\end{equation}

Note that $\tilde{\lambda}$ has the following definition.

\begin{equation}
    \tilde{L}\boldsymbol{=}\tilde{S_{0}}e^{-\tilde{S_{0}}}
\end{equation}

\begin{equation}
    \tilde{\lambda}\boldsymbol{=}-W_{0}(-\tilde{L})-1
\end{equation}

Of course, keeping in mind that the initial condition of $V(0)\boldsymbol{=}\tilde{S_{0}}+\tilde{I_{0}}$ gives the initial condition $H_{1}(0)\boldsymbol{=} 1$, equation (\ref{127}) can be easily solved for the following result.

\begin{equation}
    H_{1}(y)\boldsymbol{=}\frac{\tilde{S_{0}}-(1-y)^{\Theta}}{\tilde{S_{0}}-1}
\end{equation}

\begin{equation}
    \Theta\boldsymbol{=}\frac{\tilde{S_{0}}-1}{\tilde{\lambda}}
\end{equation}

Due to the properties of the Lambert function, $\tilde{\lambda}\boldsymbol{=}\tilde{S_{0}}-1$ when $0<\tilde{S_{0}}<1$ thereby giving $\Theta\boldsymbol{=} 1$. So, in this restricted domain, the function $H_{1}$ simplifies to the following.

\begin{equation}
    H_{1}(y)\boldsymbol{=} 1+\frac{y}{\tilde{S_{0}}-1},\,\,\,\,\,\,\,\,\,\, 0<\tilde{S_{0}}<1
\end{equation}

However, note that when $\tilde{S_{0}}>1$, $\tilde{\lambda}$ is no longer linearly related to $\tilde{S_{0}}$, and so, a power expansion of $1-y$ occurs within the expression of $H_{1}$ thereby resulting in a radius of $1$. In general, when $\tilde{S_{0}}>1$, $\Theta$ will be negative. On the downside though, $H_{1}$ will become arbitrarily large as $y\rightarrow 1$.

A further asymptotic analysis can be done to obtain $H_{2}$ while continuing to assume $\tilde{I_{0}}<<1$. This results in the limiting differential equation below.

\begin{equation}
    \frac{dH_{2}}{dy}-\left(\frac{\Theta}{y-1}\right)H_{2}\boldsymbol{=}-\left(\frac{\tilde{\lambda}+1}{\tilde{\lambda}^{3}}\right)(1-y)^{\Theta-1}-\frac{\tilde{S_{0}}(1-y)^{-1}}{2\tilde{\lambda}(\tilde{S_{0}}-1)^{2}}\left(1-(1-y)^{\Theta}\right)^2
\end{equation}

Solving the equation with the initial condition $H_{2}\boldsymbol{=} 0$ since there are no $\tilde{I_{0}}^{2}$ terms when $y\boldsymbol{=} 0$, we are now able to obtain the following solution.

\begin{equation}
    H_{2}(y)\boldsymbol{=}\left(\frac{\tilde{\lambda}+1}{\tilde{\lambda}^{3}}-\frac{\tilde{S_{0}}}{\tilde{\lambda}(\tilde{S_{0}}-1)^{2}}\right)(1-y)^{\Theta}\ln(1-y)-\frac{\tilde{S_{0}}}{2(\tilde{S_{0}}-1)^{3}}\left(1-(1-y)^{2\Theta}\right)
    \label{134}
\end{equation}

If $\tilde{S_{0}}<1$, then equation (\ref{134}) simplifies to the following well behaved form.

\begin{equation}
    H_{2}(y)\boldsymbol{=}-\frac{\tilde{S_{0}}}{2(\tilde{S_{0}}-1)^{3}}\left(1-(1-y)^{2}\right),\,\,\,\,\,\,\,\,\,\, 0<\tilde{S_{0}}<1
\end{equation}

Consequently, this now answers why the off-axis behavior of the $V$ power series even with arbitrarily small $\tilde{I_{0}}$ will still have a finite radius. In fact, it can be seen from the binomial series expansion that $\tilde{I_{0}}H_{1}$ will as expected give the correct values of the $A_{n}$ coefficients of $V$ to within the first order of $\tilde{I_{0}}$.

However, there is the other issue of what happens when $\tilde{S_{0}}>1$ since $\Theta<0$ in this regime, and so, in this scenario, both $H_{1}$ and $H_{2}$ diverge to $\pm\infty$ as $y\rightarrow 1$. Clearly, some issue has occurred in the asymptotic analysis that causes it to lose validity when $\tilde{S_{0}}>1$. Possibly,  this is related to the Taylor series of $W_{0}(x)$ having a radius of convergence equal to $1/e$ with the Lambert function also not being injective, and so, it seems likely that a different asymptotic approach can yield the correct $H_{n}$ functions in the $\tilde{S_{0}}>1$ regime.

Unfortunately though, even with an additional asymptotic analysis, neither this or higher order contributions such as $H_{2}$ alone are likely to tell us how a critical value of $\tilde{S_{0}}$ occurs as $\tilde{I_{0}}\rightarrow 0$. Most likely, a holistic analysis of the sequence of higher order terms interacting together will be needed to answer this question.

\subsection{Appendix K: Asymptotic Analysis When $\tilde{S_{0}}\boldsymbol{=} 1$}

Unfortunately, if $\tilde{S_{0}}\boldsymbol{=} 1$, then the previously discussed asymptotic analysis will fail. With this in mind, an alternative asymptotic approach based on powers of $\tilde{I_{0}}/\lambda$ will be attempted for the specific case $\tilde{S_{0}}\boldsymbol{=} 1$. This is a legitimate approach since the following limit is true in that case.

\begin{equation}
   \lim_{\tilde{I_0} \to 0} \frac{\tilde{I_{0}}}{\lambda} \boldsymbol{=} \lim_{\tilde{I_0} \to 0} \frac{1}{\partial\lambda/\partial\tilde{I_{0}}}\boldsymbol{=}\lim_{\tilde{I_0} \to 0} \frac{1}{-L W_{0}'(-L)}\boldsymbol{=}-\frac{1+W_{0}(-e^{-1})}{W_{0}(-e^{-1})}\boldsymbol{=} 0
\end{equation}

However, the approach ultimately gets modified to instead more generally involve factors of $\tilde{I_{0}}^{n}/\lambda^{m}$. This results in the following expansion.

\begin{equation}
    V\boldsymbol{=}J_{0,0}+\tilde{I_{0}}J_{1,0}+\frac{\tilde{I_{0}}}{\lambda}J_{1,1}+\frac{\tilde{I_{0}}^{2}}{\lambda^{3}}J_{2,3}+\frac{I_{0}^{3}}{\lambda^{4}}J_{3,4}+\frac{I_{0}^{3}}{\lambda^{5}}J_{3,5}+...
\end{equation}

The first two conditions are needed as follows. This leads to the following asymptotic analysis results.

\begin{equation}
    J_{0,0}\boldsymbol{=} S_{0}
\end{equation}

\begin{equation}
    J_{1,0}\boldsymbol{=} 1
\end{equation}

This leads to the limiting asymptotic component equations.

\begin{equation}
    \frac{dJ_{1,1}}{dy}\boldsymbol{=} \frac{1}{1-y}
\end{equation}

\begin{equation}
    \frac{dJ_{2,3}}{dy}\boldsymbol{=}-\frac{J_{1,1}}{2(1-y)}
\end{equation}

\begin{equation}
    \frac{dJ_{3,4}}{dy}\boldsymbol{=}-\frac{1}{6}J_{1,1}^{3}
\end{equation}

\begin{equation}
    \frac{dJ_{3,5}}{dy}\boldsymbol{=}-J_{1,1}J_{2,3}
\end{equation}

Solving for these asymptotic components gives the following results.

\begin{equation}
    J_{1,1}(y)\boldsymbol{=}-\ln(1-y)
\end{equation}

\begin{equation}
    J_{2,3}(y)\boldsymbol{=}\frac{1}{6}\ln^{3}(1-y)
\end{equation}

\begin{equation}
    J_{3,4}(y)\boldsymbol{=}-\frac{1}{24}\ln^{4}(1-y)
\end{equation}

\begin{equation}
    J_{3,5}(y)\boldsymbol{=}-\frac{1}{30}\ln^{5}(1-y)
\end{equation}

All results involve expressions proportional to $\ln^{m}(1-y)$, and consequently, since $\ln(1-y)\boldsymbol{=}\lambda T$, this expansion being computed really corresponds to the special case $\tilde{S_{0}}\boldsymbol{=} 1$ of the original power series of $T$. However, as a practical matter, this is okay. As originally noted in the contour plot of the radius of convergence of this power series, it was seen that it appeared to start doing well when $\tilde{S_{0}}\boldsymbol{=}1$ with the radius probably becoming arbitrarily large in the regime of $\tilde{I_{0}}<<1$.

\subsection{Appendix L: An Asymptotic Analysis of Small $\tilde{S_{0}}$}

Given the issues that exist when $\tilde{S_{0}}>1$ in the previous analysis, let us try an altogether different idea. Now, instead, we will consider the asymptotic regime of expansions of powers of $\tilde{S_{0}}$ in which $\tilde{S_{0}}<<$ is initially assumed to be true. In other words, we will be considering the following.

\begin{equation}
    V\boldsymbol{=} P_{0}(y)+\tilde{S_{0}}P_{1}(y)+\tilde{S_{0}}^{2}P_{2}(y)+...
\end{equation}

Taking equation (\ref{44}) to define the zeroth order contribution to this idea, we have the following while using $\lambda\boldsymbol{=}-1$ in the $\tilde{S_{0}}\rightarrow 0$ limit.

\begin{equation}
    P_{0}\boldsymbol{=}\tilde{I_{0}}e^{-T}\boldsymbol{=}\tilde{I_{0}}e^{\lambda T}\boldsymbol{=}\tilde{I_{0}}(1-y)
\end{equation}

Taking into account the zeroth contribution, we are then able to arrive at the equation of $P_{1}$ in the limit of $\tilde{S_{0}}\rightarrow 0$.

\begin{equation}
    \frac{dP_{1}}{dy}-\frac{P_{1}}{y-1}\boldsymbol{=}\frac{e^{-\tilde{I_{0}}y}}{1-y}-\tilde{I_{0}}e^{-\tilde{I_{0}}}
\end{equation}

While noting that the initial condition of $V$ is $\tilde{I_{0}}+\tilde{S_{0}}$ implies there is the initial condition $P_1(0)\boldsymbol{=} 1$, we are able to solve the above equation for $P_{1}$ giving us the following result.

\begin{equation}
    P_{1}(y)\boldsymbol{=} 1-y+\tilde{I_{0}}e^{-\tilde{I_{0}}}(1-y)\ln(1-y)+(1-y)\int_{0}^{y}\frac{e^{-\tilde{I_{0}}Y}}{(1-Y)^{2}}dY
\end{equation}

Since the integral inside $P_{1}$ is divergent in the limit of $y\rightarrow 1$, $P_{1}(1)$ can be found in the following way using L'Hospital's Rule.
\begin{equation}
    P_{1}(1)\boldsymbol{=}\lim_{y\rightarrow 1}\frac{1}{(1-y)^{-1}}\int_{0}^{y}\frac{e^{-\tilde{I_{0}}Y}}{(1-Y)^{2}}dY\boldsymbol{=}\lim_{y\rightarrow 1}\frac{1}{(1-y)^{-2}}\frac{e^{-\tilde{I_{0}}y}}{(1-y)^{2}}\boldsymbol{=} e^{-\tilde{I_{0}}}
\end{equation}

\subsection{Appendix M: Initial Thoughts and Observations on the Boundary of Convergence}

If we were to consider our asymptotic series results at $y\boldsymbol{=} 1$, we could then investigate the conditions that allow full convergence of the series constructed either in powers of $\tilde{I_{0}}$ or $\tilde{S_{0}}$. So, there is a very good chance this can allow us to directly get some understanding of why the Hershey Kiss region has its particular boundary of convergence in the parameter space.

Note that this will require the individual asymptotic components to be finite at $y=1$, and unfortunately, when $\tilde{S_{0}}>1$, $\Theta<0$ which causes the $H_{j}$ components to diverge at $y=1$ in the asymptotic power series of $\tilde{I_{0}}$. So, that domain must be explored by some other means.

Consequently, the investigation of that series must be restricted to $0<\tilde{S_{0}}<1$ since this gives us simply $\Theta\boldsymbol{=} 1$ making the $\tilde{I_{0}}$ series properly convergent. If were to employ the root test in this regime on $H_{j}$ as coefficients of $\tilde{I_{0}}^{j}$, this should then start giving us a rough estimate of the maximum $\tilde{I_{0}}$ for a given $\tilde{S_{0}}$. Of course, in the limit of $j\rightarrow\infty$, we should expect this to become arbitrarily accurate giving us the convergence boundary when $0<\tilde{S_{0}}<1$. Let us take note of the following calculations.

\begin{equation}
    |H_{1}(1)|^{-1}\boldsymbol{=}\frac{1}{\tilde{S_{0}}}-1,\,\,\,\,\,\,\,\,\,\, 0<\tilde{S_{0}}<1
\end{equation}

\begin{equation}
    |H_{2}(1)^{-1/2}|\boldsymbol{=}\sqrt{2}(1-\tilde{S_{0}})^{3/2}\tilde{S_{0}}^{-1/2},\,\,\,\,\,\,\,\,\,\, 0<\tilde{S_{0}}<1
\end{equation}

\begin{equation}
    \tilde{I}_{0,max}\boldsymbol{=}\lim_{j \to \infty} |H_{j}(1)|^{-1/j},\,\,\,\,\,\,\,\,\,\, 0<\tilde{S_{0}}<1
\end{equation}

Note that the first two of the above expressions contain singularities at $\tilde{S_{0}}\boldsymbol{=}0$ thereby most likely giving a vertical asymptote for the $\tilde{I}_{0,max}$ which is exactly what we observe happening in the upper left area of the Hershey Kiss boundary of convergence.

Since the $\tilde{I_{0}}$ series fails when $\tilde{S_{0}}>1$, the next step is to try this strategy instead on the $\tilde{S_{0}}$ series. With the first coefficient $P_{1}(y)$ being an integral function (giving $P_{1}(1)\boldsymbol{=}0$), it will require much more effort to efficiently obtain numerical values of $P_{2}(1)$, $P_{3}(1)$, and beyond. It appears likely this approach should enable us to obtain a semi-analytical justification of the existence of the critical number $\tilde{S}_{0,max}$ when $\tilde{I_{0}}\rightarrow 0$. Our expectation is that the following root test strategy will give us the critical number in question.

\begin{equation}
    \tilde{S}_{0,c}\boldsymbol{=}\lim_{\tilde{I_{0}} \to 0} \lim_{j \to \infty} |P_{j}(1)|^{-1/j}
\end{equation}

\subsection{Appendix N: A Toy Model Explaining the Behavior of the Coefficients}

In the semi-log plots of the absolute value of the coefficients (of $V$) for various scenarios, we see one of two outcomes occur. In the first scenario, the coefficients eventually if not immediately begin to grow larger without bound while also having oscillations. Showing up in the semi-log representation, we see the oscillations occur about a line with a positive slope. This scenario always occurs inside vast divergent region where $\rho<1$ occurs.

In the second scenario, the coefficients grow smaller in a damped oscillation, and in the semi-log representation, we see them oscillating about a line with a negative slope. However, at some point, after many oscillations, there is a sudden phase transition in which we suddenly have a nearly horizontal line referred to previously as a "tail". This has always been observed in the convergent region otherwise known as the Hershey Kiss area.

Based on a useful discussion with another PhD student Cade Reinburger, much of this mathematical phenomena is relatively generic to various nonlinear differential equations. In his argument, the oscillatory behavior occurs due to a dominant pair of conjugate complex singularities that determine the radius of convergence by their locations being closest to the origin in the complex plane of $y$.

However, it is the author's argument that the "tail" itself is a result of having a fixed singularity at $y\boldsymbol{=} 1$ due to the $(1-y)^{-1}$ factor present in our differential equation setup. This is what is different in our system. So, $y\boldsymbol{=} 1$ and the most important pair of complex conjugate singularities will be included in this part of the discussion.

As a simple model of our solution, the following idea of an approximation of the coefficients will be investigated. As Reinburger noted, the following identity with two simple conjugage singularities is true. Interestingly, this results in a series sum of Chebyshev functions of the Second Kind.

\begin{equation}
    \frac{1}{(\rho e^{i\phi}-y)(\rho e^{-i\phi}-y)}\boldsymbol{=}\frac{1}{\rho^{2}}\left(1+\frac{y\sin(2\phi)}{\rho\sin(\phi)}+\frac{y^{2}\sin(3\phi)}{\rho^{2}\sin(\phi)}+...\right)\boldsymbol{=}\frac{1}{\rho^{2}}\sum_{n\boldsymbol{=} 0}^{\infty}\frac{y^{n}\sin((n+1)\phi)}{\rho^{n}\sin(\phi)}
\end{equation}

Generally, even if the $1$ is the dominant radius from the $y\boldsymbol{=} 1$ singularity, a larger subdominant radius from conjugate singularities is clearly playing big role for many oscillations. So, the following toy model will be assumed.
\begin{equation}
    V\approx \frac{M}{1-y}+\frac{N}{(\rho e^{i\phi}-y)(\rho e^{-i\phi}-y)},\,\,\,\,\,\,\,\,\,\, N>>M
\end{equation}

\begin{equation}
    V\approx \sum_{n\boldsymbol{=} 0}^{\infty}\left(M+\frac{N\sin((n+1)\phi)}{\rho^{n+2}\sin(\phi)}\right)y^{n}
\end{equation}

\begin{equation}
    A_{n}\approx M+\frac{N\sin((n+1)\phi)}{\rho^{n+2}\sin(\phi)}
\end{equation}

\begin{equation}
    \log|A_{n}|\approx \log\left|M+\frac{N\sin((n+1)\phi)}{\rho^{n+2}\sin(\phi)}\right|
    \label{lasteqn}
\end{equation}

A couple of simple graphical examples of equation (\ref{lasteqn}) can be seen (Figures 20 and 21) which reproduce some of the phenomenological behavior seen from points in the convergent (see Figures 9, 12, 17, and 19) and divergent regions (see Figure 18).
\begin{figure}[h!]
    \centering
    \includegraphics[width=1\linewidth]{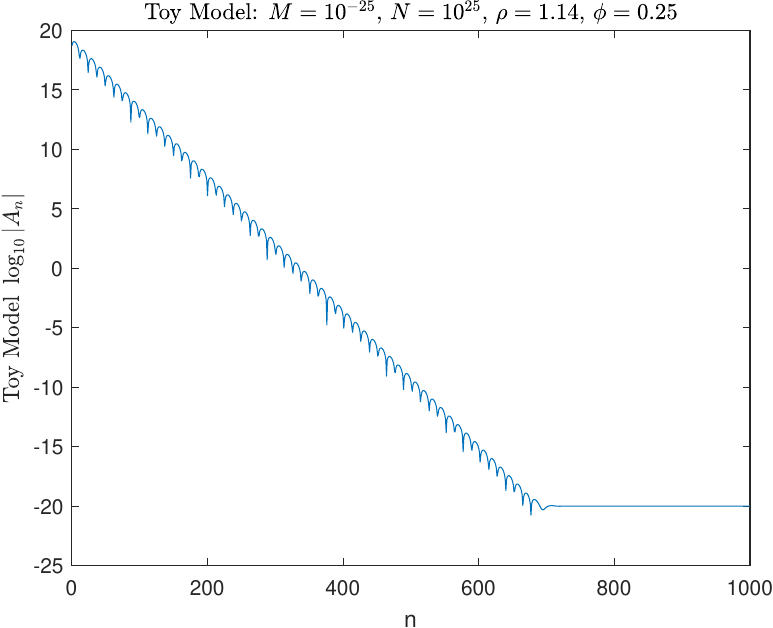}
    \caption{When the radius of convergence $\rho>1$, the Toy Model correctly reproduces a tail. In this specific case, $\rho\boldsymbol{=} 1.14$ is used.}
    \label{fig: }
\end{figure}

However, unlike Figure 17, note that Figures 8, 11, 14, and 16 actually show slightly negative slopes in their tails. This might indicate that contrary to the assumption of equation (\ref{lasteqn}), $y\boldsymbol{=} 1$ is not really a singularity and somehow factors out, but nonetheless, it would seem other singularities which are very close to $y\boldsymbol{=} 1$ (with $\rho>1$) somehow still show up. This modification of the assumption of equation (\ref{lasteqn}) would reproduce the slightly negative slopes of the tails.
\begin{figure}[h!]
    \centering
    \includegraphics[width=1\linewidth]{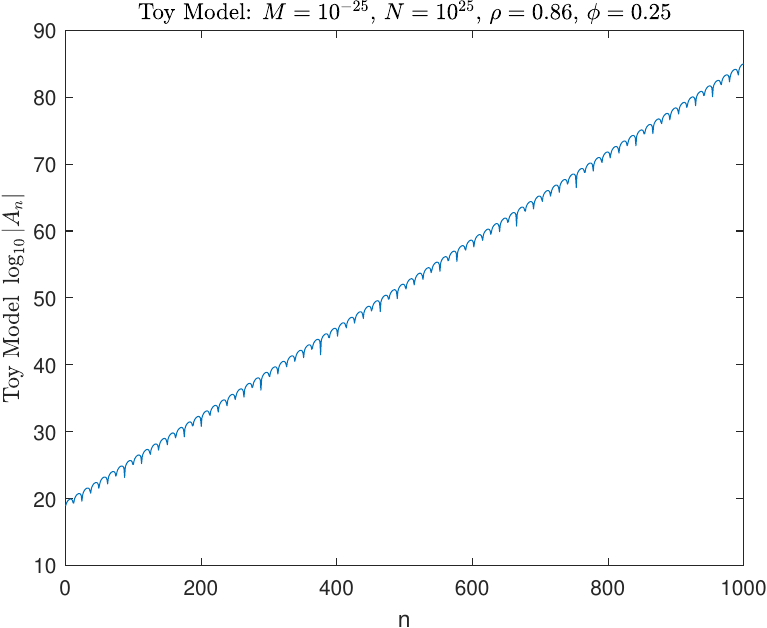}
    \caption{When the radius of convergence $\rho<1$, the Toy Model correctly reproduces divergent behavior with exponentially increasing coefficients. In this specific case, $\rho\boldsymbol{=} 0.86$ is used.}
    \label{fig: Figure 18}
\end{figure}

Finally, Figure 18 reproduces the typical behavior seen in initial conditions which are outside of the HK region.


\begin{thebibliography}{1}

\bibitem{Kermack}
W.~O. Kermack and A.~G. Mc{K}endrick.
\newblock Contribution to the mathematical theory of epidemics.
\newblock {\em Proc. Roy. Soc. London A}, 115:700--721, 1927.

\bibitem{BarlowWeinstein}
N.~S. Barlow and S.~J. Weinstein.
\newblock Accurate closed-form solution of the sir epidemic model.
\newblock {\em Physica D}, 408(132540):1--4, 2020.

\bibitem{Nastaran}
W.~Cade~Reinberger Nastaran~Naghshineh and Steven J.~Weinstein Nathaniel
  S.~Barlow, Mohamed A.~Samaha.
\newblock On the use of asymptotically motivated gauge functions to obtain
  convergent series solutions to nonlinear odes.
\newblock {\em IMA Journal of Applied Mathematics}, 88:43--66, 2023.

\bibitem{Gibbons}
A.~Gibbons.
\newblock A program for the automatic integration of differential equations
  using the method of taylor series.
\newblock {\em SIAM Rev.}, 2013.

\bibitem{Khan}
H.~Khan, R.~N. Mohapatra, K.~Vajravelu, and S.~J. Liao.
\newblock The explicit series solution of sir and sis epidemic models.
\newblock {\em Appl. Math. Comp.}, 215:653--669, 2009.

\bibitem{Rachah}
A.~Rachah and D.~F.~M. Torres.
\newblock Predicting and controlling the ebola infection.
\newblock {\em Math Method Appl. Sci.}, 40(17):6155--6164, 2017.

\bibitem{covid}
{John Hopkins University CSSE}.
\newblock Novel coronavirus (covid-19) cases.
\newblock https://github.com/CSSEGISandData/COVID-19.

\end{thebibliography}
\end{document}